\documentclass[11pt]{article}
\usepackage{a4}
\usepackage{amsmath,amssymb,amsthm}
\numberwithin{equation}{section}

\begin{document}

Acta Applicandae Math., 79, (2003), 281-326.

%\documentstyle[12pt]{article}
%\textheight=21cm
%\addtolength\oddsidemargin{-2cm}%
%\begin{document}

\title{
\Large{ Optimal with Respect to Accuracy Algorithms for
Calculation of Multidimensional Weakly Singular Integrals and
Applications to Calculation of Capacitances of Conductors of Arbitrary
Shapes
}
      }
\author{\normalsize{Ilya V. Boikov ${}^*$ and Alexander G. Ramm 
${}^{**}$}}
\date{}
\maketitle

MSC(2000) 65D32, 78A30, 78M25


{\bf Abstract:} Cubature formulas, asymptotically optimal with respect to 
accuracy, are derived 
for calculating multidimensional weakly singular
integrals. They are used for developing
a universal code for calculating capacitances of conductors
of arbitrary shapes.

{\bf 1. Introduction}

Optimal with respect to accuracy
methods for calculating  singular integrals
are being actively developed presently. They represent an important
field of computational mathematics.
Asymptotically optimal and optimal with respect
to order ( to accuracy and to complexity) algorithms for calculating
singular integrals on closed and open contours, and
multidimensional singular integrals have been constructed in [1-3]
on H\"older and Sobolev classes  of functions.

In constructing optimal with respect to accuracy methods for calculating
one-dimensional, bisingular and multidimensional singular integrals,
a general method, proposed in monograph [1], was used.
This method can be applied not only to singular integrals but also
to other integrals with moving singularities.

----------

${}^* $ Department of Higher and Applied Mathematics, Penza State University,
Krasnay Str., 40, Penza, 440026, Russia.\\
E-mail: boikov@diamond.stup.ac.ru

${}^{**} $ Department of  Mathematics, Kansas  State University,
Manhattan, KS 66506-2602, USA.\\
E-mail: ramm@math.ksu.edu

\newpage

This method allows one to construct several asymptotically optimal
and optimal with respect to order and
to accuracy algorithms for calculating hypersingular
integrals [4], the Poisson and Cauchy type integrals [5],
and multidimensional Cauchy type integrals.

   Although multidimensional weakly singular integrals are used
   in many applications, optimal methods for calculating
   these integrals are not developed.

   An exception is the book [1], where  asymptotically optimal
with respect
   to accuracy methods for calculating integrals of the form
$$
\int\limits_0^{2\pi} \int\limits_0^{2\pi} f(\sigma_1,\sigma_2)
\left|ctg \frac{\sigma_1-s_1}{2}\right|^{\gamma_1}
\left|ctg \frac{\sigma_2-s_2}{2}\right|^{\gamma_2}d\sigma_1 d\sigma_2,
$$
$0<\gamma_1$, $\gamma_2<1$,
were constructed on H\"older and Sobolev classes.

Thus, the development of optimal methods for calculating multidimensional
weakly singular integrals is an actual problem.
Construction of efficient cubature formulas for calculating weakly
singular integrals for calculating capacitances of conductors of arbitrary
shapes
by iterative methods proposed in [6] and [7] is very important in many
applications, for example, in wave scattering by small bodies of arbitrary
shapes and in antenna theory.
A bibliography on methods for calculating capacitances and
polarizability tensors is contained in [7].

In this paper the method proposed in [1] is generalized to
multidimensional
weakly singular integrals. As a result the analogs of the basic results
for singular integrals, obtained earlier, are obtained for weakly singular
integrals. Moreover, we study the
applications of optimal with respect to order cubature formulas
for calculating weakly singular integrals
on Lyapunov surfaces. Our results are used for constructing an
universal
code for calculating capacitances and polarizability tensors of
bodies of arbitrary shapes.

This paper consists of two parts.

In the first part of the paper optimal methods for calculating
integrals of the types:

$$
Kf \equiv \int\limits_0^{2\pi} \int\limits_0^{2\pi} \frac {f(\sigma_1,\sigma_2)d\sigma_1 d\sigma_2}
{\left(sin^2\left(\frac{\sigma_1-s_1}{2}\right)+sin^2\left(\frac{\sigma_2-s_2}{2}\right)\right)^{\lambda}},
\quad 0 \le s_1, s_2 \le 2\pi; \eqno (1.1)
$$
and
$$
Tf \equiv \int\limits_{-1}^1 \int\limits_{-1}^1 \frac {f(\tau_1,\tau_2)d\tau_1 d\tau_2}
{((\tau_1-t_1)^2+(\tau_2-t_2)^2)^{\lambda}}, \quad
-1 \le t_1, t_2 \le 1, \quad 0<\lambda<1, \eqno (1.2)
$$
are constructed on H\"older and Sobolev classes of functions.

Our results for integrals (1.1) can be
generalized to the integrals with other periodic kernels and functions.
The development of cubature formulas for integrals (1.1)
is of considerable interest because the results are applicable
to integrals with weakly singular kernels defined on closed Lyapunov
surfaces.

It will be clear from our arguments, that the results
can be generalized to $l-$dimensional integrals, $l=3,4,\cdots.$

The second part of this paper deals with the iterative methods for
calculating capacitances of conductors of arbitrary shapes.
A general
numerical method for calculating these capacitances is
developed, and
the results of numerical tests are given.

{\bf 2. Definitions of optimality.}

Various definitions of optimality of numerical methods
 and a detailed bibliography can be found in [8,9,10].
Let us recall the definitions of algorithms, optimal with respect to
accuracy,
for calculating  weakly singular integrals.

Consider the quadrature rule
$$
Tf=\sum\limits^{n_1}_{k_1=1} \sum\limits^{n_2}_{k_2=1}
\sum\limits^{\rho_1}_{l_1=0}\sum\limits^{\rho_2}_{l_2=0}
p_{k_1k_2l_1l_2}(t_1,t_2)f^{(l_1,l_2)}(x_{k_1},y_{k_2})+
$$
$$
+R_{n_1n_2}(f; p_{k_1k_2l_1l_2};x_{k_1},y_{k_2};t_1,t_2), \eqno (2.1)
$$
where coefficients $p_{k_1k_2l_1l_2}(t_1,t_2) $ and nodes
$ (x_{k_1},y_{k_2})$ are arbitrary.
Here $f^{(l_1,l_2)}(s_1,s_2)=\partial^{l_1+l_2} f(s_1,s_2)/
\partial s_1^{l_1} \partial s_2^{l_2}.$

The error of  quadrature rule (2.1) is defined as
$$
R_{n_1n_2}(f;p_{k_1k_2l_1l_2}; x_{k_1},y_{k_2})=
\sup\limits_{(t_1,t_2) \in [-1,1]^2}|R_{n_1n_2}(f;p_{k_1k_2l_1l_2};
x_{k_1},y_{k_2};t_1,t_2)|.
$$

The error of  quadrature rule (2.1) on the class $\Psi$ is defined as
$$
R_{n_1n_2}(\Psi;p_{k_1k_2l_1l_2}; x_{k_1},y_{k_2})=
\sup\limits_{f \in \Psi}R_{n_1n_2}(f,p_{k_1k_2l_1l_2};
x_{k_1},y_{k_2}).
$$

Define the functional
$$
\zeta_{n_1n_2}(\Psi)=\inf_{p_{k_1k_2l_1l_2};x_{k_1},y_{k_2}}
R_{n_1n_2}(\Psi;p_{k_1k_2l_1l_2};x_{k_1},y_{k_2}).
$$

The quadrature rule with the coefficients $p^*_{k_1k_2l_1l_2}$
and the nodes $(x_{k_1}^*,y_{k_2}^*)$
is optimal, asymptotically optimal, optimal with respect to order on the
class $\Psi$ among all quadrature rules of type (2.1) provided that:
$$
\frac{R_{n_1n_2}(\Psi; p^*_{k_1k_2l_1l_2}; x_{k_1}^*,y_{k_2}^*)}
{\zeta_{n_1n_2}(\Psi)}=1,\sim 1, \asymp 1, \quad n_1,n_2\to \infty.
$$

The symbol $ \alpha \asymp \beta$ means $A\alpha \le \beta \le
B\alpha,$ where $0<A,B< \infty.$

Consider the quadrature rule
$$
Tf=\sum\limits^{n}_{k=1}
p_{k}(t_1,t_2)f(M_{k})+
$$
$$
+R_{n}(f; p_{k};M_{k};t_1,t_2), \eqno (2.2)
$$
where coefficients $p_{k}(t_1,t_2) $ and nodes
$ (M_{k})$ are arbitrary.

The error of  quadrature rule (2.2) is defined as
$$
R_{n}(f;p_{k}; M_{k})=
\sup\limits_{(t_1,t_2) \in [-1,1]^2}|R_{n}(f;p_{k};M_{k};t_1,t_2)|.
$$

The error of  quadrature rule (2.2) on the class $\Psi$ is defined as
$$
R_{n}(\Psi;p_{k}; M_{k})=
\sup\limits_{f \in \Psi}R_{n}(f,p_{k};M_{k}).
$$

Define the functional
$$
\zeta_{n}(\Psi)=\inf_{p_{k};M_{k}}
R_{n}(\Psi;p_{k};M_{k}).
$$

The quadrature rule with the coefficients $p^*_{k}$
and the nodes $(M_{k}^*)$
is optimal, asymptotically optimal, optimal with respect to order on the
class $\Psi$ among all quadrature rules of type (2.2) provided that:
$$
\frac{R_{n}(\Psi; p^*_{k}; M_{k}^*)}
{\zeta_{n}(\Psi)}=1,\sim 1, \asymp 1, \quad n\to \infty.
$$

By $R_{n_1n_2}(\Psi)$ the error of optimal cubature formulas on the class
 $\Psi$ is defined. It is obvious that $R_{n_1n_2}(\Psi)=
\zeta_{n_1n_2}(\Psi).$

\hskip 50 pt {\bf 3. Classes of functions}
\vskip 10 pt

In this section, we list several classes of functions which are
used below. Some definitions are from [11,12].

A function $f$ is defined on A=[a,b] or on  A=K, where $K$ is a unit 
circle,
satisfies  the H\"older condition with  constant $M$ and exponent
$\alpha$,
or belongs to the class $H_\alpha(M), M>0,\,\, 0< \alpha \le 1,$\, if
$|f(x')-f(x'')|\le M|x'-x''|^\alpha$ \  for any $x',x'' \in A.$

Class $H_\omega,$ where $\omega(h)$ is a modulus of continuity, consists of all
functions $f\in C(A)$ with the property $|f(x_1)-f(x_2)|\le M\omega(|x_1-x_2|),
x_1,x_2 \in A.$

Class $W^r(M)$ consists of functions  $f\in C(A)$ which have continuous
derivatives
$f',f'',\dots, f^{(r-1)}$ on $A,$ and a  piecewise-continuous derivative 
$f^{(r)}$ on $A$
satisfying $max_{x\in [a,b]}|f^{(r)}(x)|\le M.$

Class  $W^r_p(M), \ r=1,2\dots, 1 \le p \le \infty,$
consists of functions $f(t),$ defined on a segment [a,b] or on $A=K,$
that have continuous derivatives
$f',f'',\dots, f^{(r-1)},$ and an integrable  derivative $f^{(r)}$
such that
$$
\left[ \int\limits_{A} |f^{(r)}(x)|^{p}dx\right]^{1/p} \le M.
$$

Class  $W^r_\alpha (M), \ r=1,2\dots, 0 <\alpha \le 1,$
consists of functions $f(t),$ defined on a segment [a,b] or on $A=K,$
which have continuous derivatives
$f',f'',\dots, f^{(r)},$
such that
$$
 |f^{(r)}(x_1)-f^{(r)}(x_2)| \le M|x_1-x_2|^\alpha.
$$

A function $f(x_1,x_2,\dots,x_l), l=2,3,\dots,$ defined on
$A=[a_1,b_2;a_2,b_2; \\ \dots;a_l,b_l]$ or on $A=K_1\times 
K_2\times\dots\times K_l,$
where $K_i, =1,2,\dots,l,$ are unit circles, satisfying H\"older
conditions with constant $M$ and exponent $\alpha_i, i=1,2,\dots,l,$
or belongs to the class $H_{\alpha_1,\dots,\alpha_l}(M), M>0, 0< \alpha\le 
1,
i=1,2,\dots,l,$ if
$$
|f(x_1,x_2,\dots,x_l)-f(y_1,y_2,\dots,y_l)|\le
M(|x_1-y_1|^{\alpha_1}+\dots+|x_l-y_l|^{\alpha_l}).
$$

Let $\omega, \omega_i,$ where $ i=1,2,\dots,l, l=1,2,\dots,$ be moduli
of continuity.

Class   $H_{\omega_1,\dots,\omega_l}(M) ,$
consists of all functions $f\in C(A),
A=[a_1,b_2;a_2,b_2;\dots;a_l,b_l]$ or $A=K_1\times K_2\times\dots\times K_l,$
with the property
$$
|f(x_1,x_2,\dots,x_l)-f(y_1,y_2,\dots,y_l)|\le
M(\omega_1(|x_1-y_1|)+\dots+\omega_l(|x_l-y_l|)).
$$

Let $H_j^\omega(A), j=1,2,3,
A=[a_1,b_2;a_2,b_2;\dots;a_l,b_l]$ or $A=K_1\times K_2\times\dots\times K_l,
l=2,3,\dots,$ be the class of functions $f(x_1,x_2,\dots,x_l)$
defined on $A$ and such that
$$
|f(x)-f(y)|\le
\omega(\rho_j(x,y)), j=1,2,3,
$$
where $x=(x_1,\dots,x_l), y=(y_1,\dots,y_l), \rho_1(x,y)=
max_{1\le i\le l}(|x_i-y_i|), \rho_2(x,y)=\sum_{i=1}^l|x_i-y_i|,
\rho_3(x,y)=[\sum_{i=1}^l|x_i-y_i|^2]^{1/2}.$

Let $H_j^\alpha(A), j=1,2,3,
A=[a_1,b_2;a_2,b_2;\dots;a_l,b_l]$ or $A=K_1\times K_2\times\dots\times K_l,
l=2,3,\dots,$ be the class of functions $f(x_1,x_2,\dots,x_l)$
defined on $A$ and such that
$$
|f(x)-f(y)|\le
(\rho_j(x,y))^\alpha, j=1,2,3.
$$

More general is the class    $H_{\rho j}^\alpha(A), j=1,2,3.$
It consists of all functions $f(x)$ which can be represented
as $f(x)=\rho(x)g(x)$, where $g(x) \in H_{ j}^\alpha(A), j=1,2,3,$ and
$ \rho(x) $ is
a nonnegative weight function.

Let $Z_j^\omega(A), j=1,2,3,$
 be the class of functions $f(x_1,x_2,\dots,x_l)$
defined on $A$ and satisfying
$$
|f(x)+f(y)-2f((x+y)/2)|\le
\omega(\rho_j(x,y)/2), j=1,2,3.
$$

Let $Z_j^\alpha(A), j=1,2,3,$
 be the class of functions $f(x_1,x_2,\dots,x_l)$
defined on $A$ and satisfying
$$
|f(x)+f(y)-2f((x+y)/2)|\le
(\rho_j(x,y)/2)^\alpha, j=1,2,3.
$$

Class    $Z_{\rho j}^\alpha(A), j=1,2,3,$
 consists of all functions $f(x)$ which can be represented
as $f(x)=\rho(x)g(x)$, where $g(x) \in Z_{ j}^\alpha(A), j=1,2,3,$ and
$ \rho(x) $ is
a nonnegative weight function.

Let $W^{r_1,\dots,r_l}(M), l=1,2,\dots,$
 be  the class of functions $f(x_1,x_2,\dots,x_l)$
defined on a domain $A,$ which have continuous partial derivatives \\
$\partial ^{|v|}f(x_1,\dots,x_l)/\partial x_1^{v_1}\dots\partial
x_l^{v_l}, 0 \le |v| \le r-1, |v|=v_1+\dots+v_l, r_i\ge v_i\ge 0, i=1,2,\dots,l,
r=r_1+\dots+r_l$ and all piece-continuous derivatives of order $r,$
satisfying $\|\partial ^{r}f(x_1,\dots,x_l)/\partial x_1^{r_1}\dots\partial
x_l^{r_l}\|_C \le M$ and $\|\partial ^{r_i}f(0,\dots,0,x_i,0,\dots,0)/
\partial x_i^{r_i}\|_C \le M, \ i=1,\dots,l.$

Let $W^{r_1,\dots,r_l}_p(M), l=1,2,\dots, 1\le p \le \infty$
 be  the class of functions $f(x_1,x_2,\dots,x_l),$
defined on a domain $A=[a_1,b_1;\dots;a_l,b_l],$
with continuous partial derivatives \\
$\partial ^{|v|}f(x_1,\dots,x_l)/\partial x_1^{v_1}\dots\partial
x_l^{v_l}, 0 \le |v| \le r-1, |v|=v_1+\dots+v_l, r_i\ge v_i\ge 0, i=1,2,\dots,l,
r=r_1+\dots+r_l,$ and all derivatives of order $r,$
satisfying
$$
\|\partial ^{r}f(x_1,\dots,x_l)/\partial x_1^{r_1}\partial x_2^{r_2}
\dots\partial x_l^{r_l}\|_{L_p(A)} \le M,
$$
$$
\|\partial ^{r_1+v_2+ \cdots +v_l}f(x_1,0,\dots,0)/\partial x_1^{r_1}\partial x_2^{v_2}
\dots\partial x_l^{v_l}\|_{L_p([a_1,b_1])} \le M, |v_2|+|v_3|+\dots+
|v_l|\le r-r_1-1;
$$
$$
...\dots ...
$$
$$
\|\partial ^{v_1 + \cdots +v_{l-1}+r_l}f(0,\dots,0,x_l)/\partial x_1^{v_1}\partial x_2^{v_2}
\dots\partial x_{l-1}^{v_{l-1}}\partial x_l^{r_l}\|_{L_p([a_l,b_l])} \le M,
|v_1|+|v_2|+\dots+|v_{l-1}|\le r-r_{l-1}-1.
$$

Let 
$A=[a_1,b_2;a_2,b_2;\dots;a_l,b_l]$ or $A=K_1\times K_2\times\dots\times K_l.$
Let  $C^{r}(M) $
be the class of functions $f(x_1,x_2,\dots,x_l)$
which are defined in $A$ and which have continuous partial derivatives
of order $r.$ Partial derivatives
of order $r$ satisfy the conditions
$$
\|\frac{ \partial ^{|v|}f(x_1,\dots,x_l)}{\partial x_1^{v_1}\dots\partial
x_l^{v_l}}\|_C\le M
$$
for any $v=(v_1,\dots,v_l),$ where $v_i\ge 0, i=1,2,\dots,l$ are integer and
$\sum_{i=1}^l v_i=r.$

By $\tilde \Psi$  we denote the set of periodic functions
of the class $\Psi.$

It is known [13] that Lyapunov spheres are defined as regions bounded
by a
finite number of closed surfaces satisfying the three Lyapunov conditions:

1. At each point of the surface a tangent plane (and, therefore, a normal)
exist.

2. If $\Theta$ is the angle between the normals at the points $m_1$
and $m_2,$ and $r$ is the distance between these points, then
$$
\Theta < Ar^{\lambda}, \quad0<\lambda \le 1,
$$
where $A$ and $\lambda$ are positive numbers which do not depend on $m_1$
and $m_2$.

3. For all points of the surface, a number $d>0$ exists such that there
is exactly one point at which a straight line, parallel to the normal at
the
surface point $m$, intersects the surface inside a sphere of
radius $d$ centered at $m.$

Let $S$ be a Lyapunov sphere, and $N$ be the exterior normal to this 
sphere.
We introduce a local system of Cartesian coordinates $(\chi,\eta,\zeta),$
whose origin is located at an arbitrary point $m_0$ of $S,$ the $\zeta$
axis
is directed along the normal $N_0$ at the point $m_0,$ and the $\chi$ and
$\eta$ axes lie in the tangential plane. In a sufficiently small
neighborhood of $m_0,$ the equation of the surface $S$ in the local
coordinates $(\chi,\eta,\zeta)$ has the form
$$
\zeta=F(\chi,\eta).
$$

{\bf Definition 4.1. [13]} {\it
The surface $S$ belongs to the class $L_k(B,\alpha)$
if $F(\chi,\eta)\in W^k_{\alpha}(B),$ and the constants $B$ and $\alpha$
do not depend on the choice of the point $m_0.$
}

{\bf 4. Auxiliary statements.}

We need the following known facts from the theory of
quadrature and cubature formulas. These facts can be found, for example,
in
[11],[12],[14], [15].

{\bf Lemma 4.1.}
{\it Let $\Psi_1$ be the class of functions $ W^r_p(1),$
$1,2,\ldots, \quad 1 \le p \le \infty,$
$0 \le t \le 1,$  $f(t) \in \Psi_1,$ and the quadrature rule
$$
\int\limits^1_0 f(t)dt=\sum\limits^n_{k=1}
p_{k}f(t_k)+R_n(f)
$$
be exact on all the polynomials of order up to $p-1$, and has error
$R_n(\Psi_1)$ on the class $\Psi_1.$
Let $\Psi_2$
be the class of functions $W^r_p(1), \quad r=1,2,\ldots,
\quad 1 \le p \le \infty, \quad a \le x \le b,$
and $g(x) \in W^r_p(1).$ Then the quadrature rule
$$
\int\limits^b_a g(x)dx=(b-a)\sum\limits^n_{k=1}
p_{k}g(a+(b-a)t_k)+R_n(g)
$$
has error $R_n(\Psi_2)$ on the class of functions $\Psi_2$ and
$$
R_n(\Psi_2)=(b-a)^{r+1-1/p}R_n(\Psi_1).
$$                  }

{\bf Theorem 4.1. [11]} {\it Among quadrature formulas
$$
\int\limits^1_0 f(x)dx=\sum\limits^m_{k=1}\sum\limits^{\rho}_{l=0}
p_{kl}f^{(l)}(x_k)+R(f) \equiv L(f)+R(f)
$$
the best formula for the class $W_p^r(1) \quad (1 \le p \le \infty)$
with  $\rho=r-1$ and $r=1,2,\cdots$,  or
$\rho=r-2$ and $ r=2,4,6,\cdots$, is the unique formula defined by
the following nodes $x^*_k$ and coefficients $p^*_{kl}$:
$$x^*_k=h(2(k-1)+[R_{rq}(1)]^{1/r}), \quad k=1,2,\ldots,m,$$
$$
p^*_{kl}=(-1)^lp^*_{ml}=h^{l+1}\left\{\frac{(-1)^l}{(l+1)!}[R_{rq}(1)]^{(l+1)/r}
+\frac{1}{r!}R_{rq}^{(r-1-1)}(1)\right\},
$$
$$
(l=0,1,\ldots,\rho),\quad p^*_{k,2v}=\frac{2h^{2v+1}}{r!}R_{rq}^{(r-2v-1)}(1), \quad
\left(k=2,3,\ldots,m-1; \quad v=0,1,\ldots,\left[\frac{r-1}{2}\right]\right),
$$
$$
p^*_{k,2v+1}=0\left(k=2,3,\ldots,m-1; \quad v=0,1,\ldots,\left[\frac{r-2}{2}\right]\right),
\quad h=2^{-1}(m-1+[R_{rq}(1)]^{1/r})^{-1},
$$
and $R_{rq}(t)$ is the Chebyshev polynomial
$t^r+\sum\limits^{r-1}_{i=0}\beta_i t^i,$
deviating least from zero in the norm $L_q(-1,1)$, where
$p^{-1}+q^{-1}=1$.
Here
$$
\zeta_n[W^r_p(1)]= R_n[W^r_p(1)]=\frac{R_{rq}(1)}{2^r r!\sqrt[q]{rq+1}(m-1+
[R_{rq}(1)]^{1/r})^r}.
$$
    }

Let a function  $f(x,y)$ be given on a rectangle $D=[a,b; c,d].$
Consider the cubature formula
$$
\iint\limits_D f(x,y)dxdy=\sum\limits^m_{k=1}\sum\limits^n_{i=1}
p_{ki}f(x_k,y_i)+R_{mn}(f),
\eqno (4.1)
$$
defined by a vector $(X,Y,P)$ of a nodes $a \le x_1<x_2<\cdots<x_m \le
b,$
$c \le y_1<y_2<\cdots<y_n \le d,$ and coefficients $p_{ki}.$

{\bf Theorem 4.2} [11]. {\it Among all quadrature formulas of the form of (4.1)
the formula
$$
\iint\limits_D f(x,y)dxdy=4hq \sum\limits^m_{k=1}\sum\limits^n_{i=1}
f(a+(2k-1)h, c+(2i-1)q)+R_{mn}(f),
$$
where $h=\frac{b-a}{2m},$ $q=\frac{d-c}{2n}, $  is optimal
on the classes $H_{\omega_1,\omega_2}(D)$ and $H^{\omega}_3(D).$
In addition
$$
R_{mn}[H_{\omega_1,\omega_2}(D)]=4mn[q\int\limits_0^h \omega_1(t)dt+
h\int\limits_0^q \omega_2(t)dt];
$$
$$
R_{mn}[H^{\omega}_3(D)]=4mn \int\limits_0^q\int\limits_0^h \omega
(\sqrt{t^2+\tau^2})dt d\tau.
$$
}

Consider the cubature formulas of the form:
$$
\iint\limits_D p(x,y)f(x,y)dxdy=\sum\limits^N_{k=1} p_kf(M_k)+R(f),
\eqno (4.2)
$$
where $p(x,y)$ is a nonnegative and bounded on $D$ function,
 $p_k,$
 $M_k(M_k \in D)$  are coefficients and nodes.

{\bf Theorem 4.3} [11].
{\it Let  $p(x,y)$ be a nonnegative bounded weight function.
If  $R_N[ H^{\alpha}_{p,j}(D)]$ and $R_N[ Z^{\alpha}_{p,j}(D)],$
where
$j=1,2,3,$ and
$0<\alpha \le 1,$ are
the  errors of optimal formulas as (4.2) on the classes
 $H^{\alpha}_{p,j}(D)$ and $ Z^{\alpha}_{p,j}(D),$ respectively,
then
$$
\lim\limits_{N\to \infty} N^{\alpha/2} R_N[H^{\alpha}_{p,j}(D)]=
\lim\limits_{N\to \infty} N^{\alpha/2} 2R_N[Z^{\alpha}_{p,j}(D)]=
$$
$$
=D_j\left[\int\int\limits_D
(p(x,y))^{2/(2+\alpha)}dxdy\right]^{(2+\alpha)/\alpha},
\quad j=1,2,3,
$$
 where
$D_1=\frac{12}{2+\alpha}\left(\frac{1}{2\sqrt{3}}\right)^{(2+\alpha)/\alpha}
\int\limits_0^{\pi/6}\frac{d\varphi}{cos^{2+\alpha}\varphi},$
$D_2=2^{1-\alpha}/(2+\alpha),$ and $D_3=2^{1-\alpha/2}/(2+\alpha).$

If  $j=2$, then the conclusion holds for $n-$dimensional cubature
formulas.}

{\bf Remark.} {\it Theorem 4.2 is generalized  to the case of
unbounded weights}
$p(x,y)$ {\it  in [2].  }

In  this  paper  we will  use the following result
(see [16]):
    \par
{\bf Lemma 4.4 }.
{\it Let $H$ be a linear metric space,
$F$ be a bounded, closed, convex, centrally symmetric set with center of  
symmetry $\theta$ at the origin, and
 $L(f),l_1(f),\dots,l_N(f),$  be some linear
functionals.
Let $S(l_1(f),\dots,l_N(f))$
be some method for calculating the functional $L(f)$ using
functionals $(l_1(f),\dots,l_N(f)),$ and $\mathcal {S}$ be the set of all
such methods.
Then the numbers $D_1,\dots,D_N$ exist such that
$$\sup_{f\in F}|L(f)-\sum\limits_{k=1}^ND_kl_k(f)|=
\inf_{\mathcal {S}}\sup_{f\in F}|L(f)-S(l_1(f),\dots,l_N(f))|.
\eqno (4.3)
$$
This means that among the best methods for calculating
functional $L(f)$:
$$ L(f)\approx S(l_1(f),\dots,l_N(f)), \eqno (4.4)$$
there is a linear method.}\par

{\bf Proof.}
Let us associate with each $f \in F$  a point
$(L(f), l_1(f),\ldots,l_N(f)).$ Let $Y$ be a set of all such
points $(y_0,\ldots,y_N)$ for $f \in F.$

From our assumptions, it follows that $Y$
is a closed centrally symmetric set with the center of symmetry
at the origin.

Let $(y_0,0,\ldots,0)$ be an extremal point of the set $Y$, and
$$
D_0=\sup_{(z,0,\ldots,0) \in Y} z =y_0.
$$
Because $F$ is bounded, one has $D_0< \infty,$ 
and because $F$ is convex and centrally symmetric with respect to
the origin, one has $D_0>0$. 

Draw the support plane for the set $Y$ through the point
$(D_0, 0,\ldots,0):$
$$
(y_0-D_0)+\sum\limits^N_{j=1}C_j y_j=0.
$$
Since $Y$ is centrally symmetric with respect to the
origin, the plane
$$
(y_0+D_0)+\sum\limits^N_{j=1}C_j y_j=0
$$
is also a support plane for $Y$, and $Y$ lies between
these two planes.

Hence, we have for the points of $Y$ the inequality:
$$
|y_0-\sum\limits^N_{j=1} D_j y_i| \le D_0, \quad D_j=-C_j.
$$

The definition of $y_i$ implies
$$
\sup_{f\in F}|L(f)-\sum\limits^N_{j=1} D_j l_j(f)| \le D_0.    \eqno (4.5)
$$

Let $f_0$ be an element $F$ corresponding the point
$(D_0, 0,\ldots,0).$ Then
$S(l_1( \pm f_0),\ldots,l_N( \pm f_0))=S(0,\ldots,0).$
%$L(f_0)-L(-f_0)=2D_0.$ Therefore either for $f_0$ or
%for $-f_0$ the error in (4.3) is not less than $D_0.$ 
The right-hand side of (4.3) is not less than 
$$\inf_{\mathcal 
{S}}\max |L(f_0)-S(0,....,0)|=\inf_a\max \{|D_0-a|,|D_0+a|\}=D_0.$$
This and
(4.5) imply that the right-hand side in (4.3)  is not less
that the left-hand one. But the right-hand side of (4.3)
 can not be more than the left-hand side of (4.3) because a
set of methods $\mathcal {S}$  contains linear methods. Lemma 4.4 is 
proved.
 $\blacksquare$

{\bf Corollary.} {\it Among all functions for which the optimal
method for calculating $L(t)$ has the greatest error
for a given set
of functionals,  there exists
a function satisfying the conditions $l_1(f)= \cdots =l_N(f)=0.$
}
It follows
from the proof that such a function is the function $f_0.$

{\bf 5. Optimal methods for calculating integrals of the form  (1.1).}

{\bf 5.1. Lower bounds for the functionals $\zeta_{nm}$ and $\zeta_N$.}

In this Section we derive lower bounds for the functionals $\zeta_{nm}$
and $\zeta_N$,
defined in Section 2 , for calculating integrals
(1.1) by the cubature formulas
$$
Kf=\sum\limits^{n_1}_{k_1=1} \sum\limits^{n_2}_{k_2=1}
\sum\limits^{\rho_1}_{l_1=0} \sum\limits^{\rho_2}_{l_2=0}p_{k_1k_2l_1l_2}
(s_1,s_2)f^{(l_1,l_2)}(x_{k_1},x_{k_2})+
$$
$$
+R_{n_1n_2}(f;p_{k_1k_2l_1l_2};x_{k_1},x_{k_2};s_1,s_2), \eqno (5.1)
$$
and
$$
Kf=\sum\limits^N_{k=1}p_k(s_1,s_2)f(M_k)+R_N(f;p_k; M_k;s_1,s_2)
\eqno (5.2)
$$
on H\"older and Sobolev classes.

{\bf Theorem  5.1.} {\it Let
$\Psi=H_{\omega_1,\omega_2}(D)$ or $\Psi=H^{\omega}_3
(D),$ and calculate integral (1.1) by formula
(5.1) with
$\rho_1=\rho_2=0.$ Then the inequality
$$
\zeta_{n_1n_2}[\Psi] \ge \frac{\gamma}{\pi^2}n_1n_2[q\int\limits_0^h \omega_1(t)dt+
h\int\limits_0^q \omega_2(t)dt],
$$
where $q=\frac{\pi}{n_2},$ $h=\frac{\pi}{n_1},$ and
$$\gamma:=\int\limits_0^{2\pi}\int\limits_0^{2\pi}\frac{ds_1ds_2}{(sin^2(s_1/2)+
sin^2(s_2/2))^{\lambda}}
\eqno (5.1')
$$
is valid.              }

{\bf Corollary}.
{\it Let  $\Psi=H_{\alpha\alpha}(D)$ or $\Psi=H^{\alpha}_3(D),$ and
calculate integral (1.1) by  formula  (5.1) with
$n_1=n_2=n$ and $\rho_1=\rho_2=0.$ Then the inequality
$$
\zeta_{nn}[\Psi] \ge \frac{2\gamma \pi^{\alpha}}{(1+\alpha)n^{\alpha}}
$$
is valid.           }

{\bf Proof of Theorem  5.1.} Denote by $\psi(s_1,s_2)$
a nonnegative function belonging to the class $H_{\omega_1\omega_2}(1)$
and vanishing at the nodes $(x_{k_1},x_{k_2}),$ $1 \le k_1 \le n_1,$
$1 \le k_2 \le n_2.$

One has:
$$
R_{n_1n_2}(\psi; p_{k_1k_2}; x_{k_1},x_{k_2}) \ge
$$
$$
\ge\frac{1}{4\pi^2}
\int\limits_0^{2\pi}\int\limits_0^{2\pi}\left(
\int\limits_0^{2\pi}\int\limits_0^{2\pi}
\frac{\psi(\sigma_1,\sigma_2)d\sigma_1d\sigma_2}{[sin^2((\sigma_1-s_1)/2)+
sin^2((\sigma_2-s_2)/2))]^\lambda}\right)
ds_1ds_2=
$$
$$
=\frac{1}{4\pi^2}\int\limits_0^{2\pi}\int\limits_0^{2\pi}
\psi(\sigma_1,\sigma_2)\left(\int\limits_0^{2\pi}\int\limits_0^{2\pi}
\frac{ds_1ds_2}{[sin^2((\sigma_1-s_1)/2)+sin^2((\sigma_2-s_2)/2))]^\lambda}
\right)d\sigma_1d\sigma_2=
$$
$$
=\frac{1}{4\pi^2}\int\limits_0^{2\pi}\int\limits_0^{2\pi}
\frac{ds_1ds_2}{[sin^2 (s_1/2)+sin^2 (s_2/2)]^\lambda}
\int\limits_0^{2\pi}\int\limits_0^{2\pi}
\psi(s_1,s_2)ds_1ds_2. \eqno (5.3)
$$

From Lemma 4.4 and Theorem 4.2 one concludes that  the following
inequality
$$
R_{n_1n_2}(\psi;p_{k_1k_2};x_{k_1},x_{k_2}) \ge
\frac{\gamma}{\pi^2}n_1n_2\left [
q\int\limits_0^h \omega_1(t)dt +h\int\limits_0^q\omega_2(t)dt\right],
h=\frac{\pi}{n_1}, \quad q=\frac{\pi}{n_2}
$$
holds for arbitrary weights
 $p_{k_1k_2}$ and nodes $(x_{k_1},x_{k_2})$
and
$$
\zeta_{nn}(\Psi) \ge \frac{\gamma}{\pi^2}n_1n_2[q\int\limits_0^h
\omega_1(t)dt+
h\int\limits_0^q \omega_2(t)dt].
$$

Theorem 5.1 is proved.           $\blacksquare$

{\bf Theorem 5.2.}
{\it Let  $\Psi=H_i^{\alpha}$ or $\Psi=Z_i^{\alpha},$
$i=1,2,3,$ and calculate the integral $Kf$ by cubature
formula (5.2).
Then
$$
\zeta_N[H_i^{\alpha}]=2\zeta_N[Z_i^{\alpha}]=(1+o(1))\gamma(4\pi^2)^{2/\alpha}
D_iN^{-\alpha/2},
$$
where  $D_1=\frac{12}{2+\alpha}\left(\frac{1}{2\sqrt{3}}\right)^
{(\alpha+2)/2}\int\limits_0^{\pi/6}\frac{d\varphi}{cos^{2+\alpha}\varphi},$\,
$D_2=\frac{2}{2^{\alpha}(2+\alpha)},$ and $D_3=\frac{2^{1-\alpha/2}}
{2+\alpha}.$          }

{\bf Proof.} The proof of Theorem 5.2 is similar to the proof of Theorem
5.1, with some
 difference is in the estimation of the integral
$\int\limits_0^{2\pi}\int\limits_0^{2\pi}
\psi(s_1,s_2)ds_1ds_2,$ where the function $\psi(s_1,s_2)$
belongs to the class
$H_i^{\alpha}$ (or $Z_i^{\alpha}$), is nonnegative in the domain
$D=[0,2\pi]^2$,
and vanishes at  $N$ nodes $M_k,$ $k=1,2,\ldots,N.$

Using Lemma 4.4 and Theorem 4.3, one checks that the
inequalities
$$
\inf\limits_{M_k}\sup\limits_{\psi \in H_i^\alpha, \psi(M_k)=0}
\int\limits_0^{2\pi}\int\limits_0^{2\pi}\psi(s_1,s_2)ds_1ds_2=
(1+o(1))D_i(4\pi^2)^{(2+\alpha)/\alpha}N^{-\alpha/2},
$$
$$
\inf\limits_{M_k}\sup\limits_{\psi \in Z_i^\alpha, \psi(M_k)=0}
\int\limits_0^{2\pi}\int\limits_0^{2\pi}\psi(s_1,s_2)ds_1ds_2=
(1+o(1))\frac{1}{2}D_i(4\pi^2)^{(2+\alpha)/\alpha}N^{-\alpha/2}
$$
hold for arbitrary
$M_k \in D,$ $k=1,2,\ldots,N.$

Substituting these values into inequality (5.3), we complete the proof
of Theorem 5.2. $\blacksquare$

{\bf Theorem 5.3.}
{\it Let $\Psi=\tilde C_2^r(1)$, and calculate the integral $Kf$
by formula  (5.1) with $\rho_1=\rho_2=0,$ and $n_1=n_2=n.$ Then
$$
\zeta_{nn}[\Psi]\ge(1+o(1))\frac{2\gamma K_r}{n^r},
$$
where $K_r$ is the Favard constant.}

{\bf Proof.} Let
$$
\psi(s_1,s_2)=\psi_1(s_1)+\psi_2(s_2),
$$
where $0\leq \psi_1(s)\in W^r(1)$
vanishes at the nodes $x_k,$ $k=1,2,\ldots,n,$
and $0\leq \psi_2(s)\in W^r(1)$
vanishes at the nodes $y_k,$ $k=1,2,\ldots,n.$

According to [11], for arbitrary nodes $x_k,$ $k=1,2,\ldots,n$
one has:
$$
\int\limits_0^{2\pi}\psi_i(s)ds \ge \frac{2\pi K_r}{n^r}, \  i=1,2.
$$

Thus, the inequality
$$
\int\limits_0^{2\pi}\int\limits_0^{2\pi} \psi(s_1,s_2)ds_1ds_2 \ge
\frac{8\pi^2K_r}{n^r}
$$
holds for arbitrary nodes $(x_1,\ldots,x_n)$ and
$(y_1,\ldots,y_n).$

The conclusion  of Theorem 5.3 follows from
this inequality and from  (5.3).
$\blacksquare$

{\bf Theorem 5.4.}
{\it Let $\Psi=W_p^{r,r}(1),$ $r=1,2,\ldots,$ $1\le p \le \infty,$ and
calculate the integral $Kf$ by formula (5.1) with
$\rho_1=\rho_2=r-1$ and $n_1=n_2=n.$
Then
$$
\zeta_{nn}[\Psi] \ge (1+o(1))
\frac{2^{1/q}\pi^{r-1/p}R_{rq}(1)}{r!(rq+1)^{1/q}(n-1+[R_{rq}(1)]^{1/r})^r}
\gamma,
$$
where $R_{rq}(t)$ is a polynomial of degree $r$, least deviating
from zero in $L_q([-1,1]).$}

{\bf Proof.} Let $L=[\frac n{\log n}].$
  Take an additional set of nodes
$(\xi_k,\xi_l),$ $\xi_k=\frac{2\pi k}{L},$
$k,l=0,1,\ldots,L-1.$
By $(v_i,w_j),$ $i,j=0,1,\ldots,N-1,$ $N=n+L,$
denote the union of the sets  $(x_k,y_l)$ and $(\xi_i ,\xi_j).$
Let $\psi(s_1,s_2)=\psi_1(s_1)+\psi_2(s_2),$ where $\psi_1(s)\in W^r_p(1)$
vanishes with its derivatives up to the order  $r-1$
at the nodes  $v_i,$ $i=0,1,\ldots,N-1,$ and $\psi_2(s)\in W_p^{(r)}(1)$
vanishes  with its derivatives up to order  $r-1$
 at the nodes $w_j,$ $j=0,1,\ldots,N-1.$ Assume
$\int\limits_{v_i}^{v_{i+1}}\psi_1(s)ds>0, \,\,\,
i=0,1,\ldots,N-1,$ and
$\int\limits_{w_j}^{w_{j+1}}\psi_2(s)ds>0, \,\,\,
j=0,1,\ldots,N-1,$
where $v_{N}=2\pi$ and $w_{N}=2\pi.$

Let
$$
{\bf h(s_1,s_2,\sigma_1,\sigma_2):=}
\left\{
\begin{array}{cc}
0, \quad if \quad (\sigma_1,\sigma_2)=(s_1,s_2), \\
\frac{1}{(sin^2((\sigma_1-s_1)/2)+sin^2((\sigma_2-s_2)/2))^\lambda},
\quad otherwise,
\end{array}
\right.
$$
$$
{\bf\psi^+(s_1,s_2)=}
\left\{
\begin{array}{cc}
\psi(s_1,s_2), \quad if \quad \psi(s_1,s_2) \ge 0,  \\
0,             \quad if \quad \psi(s_1,s_2) < 0,
\end{array}
\right.
$$
$$
{\bf\psi^-(s_1,s_2)=}
\left\{
\begin{array}{cc}
0,             \quad if \quad \psi(s_1,s_2) \ge 0,\\
-\psi(s_1,s_2), \quad if \quad \psi(s_1,s_2) < 0.
\end{array}
\right.
$$

For each value  $(\xi_i,\xi_j),$ $i,j=0,1,\ldots,N-1,$
we have (with $N=N_1=N_2=L$):
$$
\int\limits_0^{2\pi}\int\limits_0^{2\pi} h(\xi_i,\xi_j,\sigma_1,\sigma_2)
\psi(\sigma_1,\sigma_2)d\sigma_1 d\sigma_2=
$$
$$
=\sum\limits_{k=0}^{N-1}\sum\limits_{l=0}^{N-1}\int\limits_{\xi_k}^{\xi_{k+1}}
\int\limits_{\xi_l}^{\xi_{l+1}} h(\xi_i,\xi_j,\sigma_1,\sigma_2)
\psi(\sigma_1,\sigma_2)d\sigma_1 d\sigma_2=
$$
$$
=\sum\limits_{k=0}^{N-1}\sum\limits_{l=0}^{N-1}\int\limits_{\xi_k}^{\xi_{k+1}}
\int\limits_{\xi_l}^{\xi_{l+1}} h(\xi_i,\xi_j,\sigma_1,\sigma_2)
\psi^+(\sigma_1,\sigma_2)d\sigma_1 d\sigma_2-
$$
$$
-\sum_{k=0}^{N-1}\sum\limits_{l=0}^{N-1}\int\limits_{\xi_k}^{\xi_{k+1}}
\int\limits_{\xi_l}^{\xi_{l+1}} h(\xi_i,\xi_j,\sigma_1,\sigma_2)
\psi^-(\sigma_1,\sigma_2)d\sigma_1 d\sigma_2 \ge
$$
$$
\ge\sum\limits_{k=i+1}^{i+[(N_1-1)/2]}\sum\limits_{l=j+1}^{j+[(N_2-1)/2]}
h(\xi_i,\xi_j, \xi_{k+1},\xi_{l+1})\int\limits_{\xi_k}^{\xi_{k+1}}
\int\limits_{\xi_l}^{\xi_{l+1}}\psi^+(\sigma_1,\sigma_2)d\sigma_1 d\sigma_2+
$$
$$
+\sum\limits_{k=i+1}^{i+[(N_1-1)/2]}\sum\limits^{j-1}_{l=j-[(N_2-1)/2]}
h(\xi_i,\xi_j, \xi_{k+1},\xi_{l-1})\int\limits_{\xi_k}^{\xi_{k+1}}
\int\limits_{\xi_{l-1}}^{\xi_l}\psi^+(\sigma_1,\sigma_2)d\sigma_1 d\sigma_2+
$$
$$
+\sum\limits^{i-1}_{k=i-[(N_1-1)/2]}\sum\limits_{l=j+1}^{j+[(N_2-1)/2]}
h(\xi_i,\xi_j, \xi_{k-1},\xi_{l+1})\int\limits^{\xi_k}_{\xi_{k-1}}
\int\limits_{\xi_l}^{\xi_{l+1}}\psi^+(\sigma_1,\sigma_2)d\sigma_1 d\sigma_2+
$$
$$
+\sum\limits^{i-1}_{k=i-[(N_1-1)/2]}\sum\limits^{j-1}_{l=j-[(N_2-1)/2]}
h(\xi_i,\xi_j, \xi_{k-1},\xi_{l-1})\int\limits^{\xi_k}_{\xi_{k-1}}
\int\limits^{\xi_l}_{\xi_{l-1}}\psi^+(\sigma_1,\sigma_2)d\sigma_1 d\sigma_2-
$$
$$
-\sum\limits_{k=i+1}^{i+[(N_1-1)/2]}\sum\limits_{l=j+1}^{j+[(N_2-1)/2]}
h(\xi_i,\xi_j, \xi_k,\xi_l)\int\limits_{\xi_k}^{\xi_{k+1}}
\int\limits_{\xi_l}^{\xi_{l+1}}\psi^-(\sigma_1,\sigma_2)d\sigma_1 d\sigma_2-
$$
$$
-\sum\limits_{k=i+1}^{i+[(N_1-1)/2]}\sum\limits^{j-1}_{l=j-[(N_2-1)/2]}
h(\xi_i,\xi_j, \xi_k,\xi_l)\int\limits_{\xi_k}^{\xi_{k+1}}
\int\limits^{\xi_l}_{\xi_{l-1}}\psi^-(\sigma_1,\sigma_2)d\sigma_1 d\sigma_2-
$$
$$
-\sum\limits^{i-1}_{k=i-[(N_1-1)/2]}\sum\limits_{l=j+1}^{j+[(N_2-1)/2]}
h(\xi_i,\xi_j, \xi_k,\xi_l)\int\limits^{\xi_k}_{\xi_{k-1}}
\int\limits_{\xi_l}^{\xi_{l+1}}\psi^-(\sigma_1,\sigma_2)d\sigma_1 d\sigma_2-
$$
$$
-\sum\limits^{i-1}_{k=i-[(N_1-1)/2]}\sum\limits^{j-1}_{l=j-[(N_2-1)/2]}
h(\xi_i,\xi_j, \xi_k,\xi_l)\int\limits^{\xi_k}_{\xi_{k-1}}
\int\limits^{\xi_l}_{\xi_{l-1}}\psi^-(\sigma_1,\sigma_2)d\sigma_1 d\sigma_2=
$$
$$
=\sum\limits_{k=i+1}^{i+[(N_1-1)/2]}\sum\limits_{l=j+1}^{j+[(N_2-1)/2]}
h(\xi_i,\xi_j, \xi_{k+1},\xi_{l+1})\int\limits_{\xi_k}^{\xi_{k+1}}
\int\limits_{\xi_l}^{\xi_{l+1}}\psi(\sigma_1,\sigma_2)d\sigma_1 d\sigma_2+
$$
$$
+\sum\limits_{k=i+1}^{i+[(N_1-1)/2]}\sum\limits^{j-1}_{l=j-[(N_2-1)/2]}
h(\xi_i,\xi_j, \xi_{k+1},\xi_{l-1})\int\limits_{\xi_k}^{\xi_{k+1}}
\int\limits^{\xi_l}_{\xi_{l-1}}\psi(\sigma_1,\sigma_2)d\sigma_1 d\sigma_2+
$$
$$
+\sum\limits^{i-1}_{k=i-[(N_1-1)/2]}\sum\limits_{l=j+1}^{j+[(N_2-1)/2]}
h(\xi_i,\xi_j, \xi_{k-1},\xi_{l+1})\int\limits^{\xi_k}_{\xi_{k-1}}
\int\limits_{\xi_l}^{\xi_{l+1}}\psi(\sigma_1,\sigma_2)d\sigma_1 d\sigma_2+
$$
$$
+\sum\limits^{i-1}_{k=i-[(N_1-1)/2]}\sum\limits^{j-1}_{l=j-[(N_2-1)/2]}
h(\xi_i,\xi_j, \xi_{k-1},\xi_{l-1})\int\limits^{\xi_k}_{\xi_{k-1}}
\int\limits^{\xi_l}_{\xi_{l-1}}\psi(\sigma_1,\sigma_2)d\sigma_1 d\sigma_2-
$$
$$
-\sum\limits_{k=i+1}^{i+[(N_1-1)/2]}\sum\limits_{l=j+1}^{j+[(N_2-1)/2]}
(h(\xi_i,\xi_j, \xi_k,\xi_l)-h(\xi_i,\xi_j, \xi_{k+1},\xi_{l+1}))\times
$$
$$
\times \int\limits_{\xi_k}^{\xi_{k+1}}\int\limits_{\xi_l}^{\xi_{l+1}}\psi^-
(\sigma_1,\sigma_2)d\sigma_1 d\sigma_2-
$$
$$
-\sum\limits_{k=i+1}^{i+[(N_1-1)/2]}\sum\limits^{j-1}_{l=j-[(N_2-1)/2]}
(h(\xi_i,\xi_j, \xi_k,\xi_l)-h(\xi_i,\xi_j, \xi_{k+1},\xi_{l-1}))\times
$$
$$
\times\int\limits^{\xi_k}_{\xi_{k+1}}\int\limits^{\xi_l}_{\xi_{l-1}}
\psi^-(\sigma_1,\sigma_2)d\sigma_1 d\sigma_2-
$$
$$
-\sum\limits^{i-1}_{k=i-[(N_1-1)/2]}\sum\limits_{l=j+1}^{j+[(N_2-1)/2]}
(h(\xi_i,\xi_j, \xi_k,\xi_l)-h(\xi_i,\xi_j, \xi_{k-1},\xi_{l+1}))\times
$$
$$
\times\int\limits^{\xi_k}_{\xi_{k-1}}\int\limits_{\xi_l}^{\xi_{l+1}}
\psi^-(\sigma_1,\sigma_2)d\sigma_1 d\sigma_2-
$$
$$
-\sum\limits^{i-1}_{k=i-[(N_1-1)/2]}\sum\limits^{l=j-1}_{l=j-[(N_2-1)/2]}
(h(\xi_i,\xi_j, \xi_k,\xi_l)-h(\xi_i,\xi_j, \xi_{k-1},\xi_{l-1}))\times
$$
$$
\times\int\limits^{\xi_k}_{\xi_{k-1}}\int\limits^{\xi_l}_{\xi_{l-1}}
\psi^-(\sigma_1,\sigma_2)d\sigma_1 d\sigma_2=
$$
$$
=J_1+J_2+J_3+J_4+I_1+I_2+I_3+I_4.
$$

Let us estimate the integral
$$
\left|\int\limits_{\xi_k}^{\xi_{k+1}}\int\limits_{\xi_l}^{\xi_{l+1}}
\psi^-(\sigma_1,\sigma_2)d\sigma_1 d\sigma_2\right| \le
\int\limits_{\xi_k}^{\xi_{k+1}}\int\limits_{\xi_l}^{\xi_{l+1}}
|\psi^-(\sigma_1,\sigma_2)|d\sigma_1 d\sigma_2 \le
$$
$$
\le\int\limits_{\xi_k}^{\xi_{k+1}}\int\limits_{\xi_l}^{\xi_{l+1}}
|\psi(\sigma_1,\sigma_2)|d\sigma_1 d\sigma_2 \le
(\xi_{l+1}-\xi_l)\int\limits_{\xi_k}^{\xi_{k+1}}|\psi_1(\sigma)|d\sigma+
(\xi_{k+1}-\xi_k)\int\limits_{\xi_l}^{\xi_{l+1}}|\psi_2(\sigma)|d\sigma \le
$$
$$
\le 2\left(\frac{2\pi }{L}\right)^{r+2}\frac{1}{r!},
$$
where we have used the fact that the functions
$\psi_1(s)$ and $\psi_2(s)$ on the segments $[\xi_k,\xi_{k+1}]$ and
$[\xi_l,\xi_{l+1}]$ vanish with derivatives up to order
$r-1$.

Now let us estimate the sum:
$$
\sum\limits_{k=i+1}^{i+[(L-1)/2]}\sum\limits_{l=j+1}^{j+[(L-1)/2]}
|h(\xi_i,\xi_j, \xi_{k+1},\xi_{l+1})-h(\xi_i,\xi_j, \xi_k,\xi_l)|=
$$
$$
=\sum\limits_{k=i+1}^{i+[(L-1)/2]}\sum\limits_{l=j+1}^{j+[(L-1)/2]}
\left|\frac{1}{\left(sin^2\frac{\pi (k+1-i)}{L}+
sin^2\frac{\pi (l+1-j)}{L}\right)^\lambda}-\right.
$$
$$
\left.-\frac{1}{\left(sin^2\frac{2\pi (k-i)}{L}+
sin^2\frac{2\pi (l-j)}{L}\right)^\lambda}\right| \le
$$
$$
\le\frac cL\sum\limits_{k=i+1}^{i+[(L-1)/2]}\sum\limits_{l=j+1}^{j+[(L-1)/2]}
\frac1{\left(sin^2\frac{\pi(k-i)}{L}+sin^2\frac{\pi(l-j)}{L}\right)^{1+
\lambda}}
\left|\frac{k-i}{L}+\frac{l-j}{L}\right|\le
$$
$$
\le\frac{c}{L}
\sum\limits_{k=i+1}^{i+[(L-1)/2]}\sum\limits_{l=j+1}^{j+[(L-1)/2]}
\frac{L^{2+2\lambda}}{\left((k-i)^2+(l-j)^2\right)^{1+\lambda}}
\frac{(k-i)+(l-j)}{L} \le
$$
$$
\le c\left(L\right)^{2\lambda}\left(\sum\limits_{l=1}^{[(L-1)/2]}
\frac{1}{l^{2\lambda}}+\sum\limits_{k=1}^{[(L-1)/2]}\frac{1}{k^{2\lambda}}
\right) \le
$$
$$
 \le c\left(L\right)^{2\lambda}
\left\{
\begin{array}{ccc}
L^{1-2\lambda} \quad if \quad \lambda<\frac{1}{2},\\
\log L           \quad if \quad \lambda=\frac{1}{2},  \\
1              \quad if \quad \lambda>\frac{1}{2}.    \\
\end{array}
\right.
$$
By $c>0$ various estimation constants are denoted.
Thus
$$
I_1=o\left(\frac{1}{n^r}\right).
$$

The expressions  $I_2,I_3,$ and $I_4$ are esimated similarly.

From the definition of the function $\psi(s_1,s_2)$ it follows
that the error of cubature formula (5.1) for $s_1=\xi_i,$ $s_2=\xi_j$
can be estimated as follows:
$$
R(\psi,   ,\xi_i,\xi_j)=\int\limits_0^{2\pi}\int\limits_0^{2\pi}
\psi(\sigma_1,\sigma_2)h(\xi_i,\xi_j,\sigma_1,\sigma_2)d\sigma_1d\sigma_2 \ge
o\left(\frac{1}{n^r}\right)+
$$
$$
+\sum\limits_{k=i+1}^{i+[(L-1)/2]}\sum\limits_{l=j+1}^{j+[(L-1)/2]}
h(\xi_i,\xi_j, \xi_{k+1},\xi_{l+1})\int\limits_{\xi_k}^{\xi_{k+1}}
\int\limits_{\xi_l}^{\xi_{l+1}}\psi(\sigma_1,\sigma_2)d\sigma_1 d\sigma_2+
$$
$$
+\sum\limits_{k=i+1}^{i+[(L-1)/2]}\sum\limits^{j-1}_{l=j-[(L-1)/2]}
h(\xi_i,\xi_j, \xi_{k+1},\xi_{l-1})\int\limits_{\xi_k}^{\xi_{k+1}}
\int\limits^{\xi_l}_{\xi_{l-1}}\psi(\sigma_1,\sigma_2)d\sigma_1 d\sigma_2+
$$
$$
+\sum\limits^{i-1}_{k=i-[(L-1)/2]}\sum\limits_{l=j+1}^{j+[(L-1)/2]}
h(\xi_i,\xi_j, \xi_{k-1},\xi_{l+1})\int\limits^{\xi_k}_{\xi_{k-1}}
\int\limits_{\xi_l}^{\xi_{l+1}}\psi(\sigma_1,\sigma_2)d\sigma_1 d\sigma_2+
$$
$$
+\sum\limits^{i-1}_{k=i-[(L-1)/2]}\sum\limits^{j-1}_{l=j-[(L-1)/2]}
h(\xi_i,\xi_j, \xi_{k-1},\xi_{l-1})\int\limits^{\xi_k}_{\xi_{k-1}}
\int\limits^{\xi_l}_{\xi_{l-1}}\psi(\sigma_1,\sigma_2)d\sigma_1 d\sigma_2.
$$

Averaging the above inequality over $i$ and $j,$ one gets:
$$
R_{nn}[\Psi] \ge \sup_{\psi \in \Psi}
\max\limits_{i,j}R_{nn}(\psi,   ,\xi_i,\xi_j) \ge
$$
$$
\ge \frac{1}{L^2}\sum\limits^{L-1}_{i=0}\sum\limits^{L-1}_{j=0}
\left[\sum\limits_{k=i+1}^{i+[(L-1)/2]}\sum\limits_{l=j+1}^{j+[(L-1)/2]}
h(\xi_i,\xi_j, \xi_{k+1},\xi_{l+1})\int\limits_{\xi_k}^{\xi_{k+1}}
\int\limits_{\xi_l}^{\xi_{l+1}}\psi(\sigma_1,\sigma_2)d\sigma_1 d\sigma_2+
\right.
$$
$$
+\sum\limits_{k=i+1}^{i+[(L-1)/2]}\sum\limits^{j-1}_{l=j-[(L-1)/2]}
h(\xi_i,\xi_j, \xi_{k+1},\xi_{l-1})\int\limits_{\xi_k}^{\xi_{k+1}}
\int\limits^{\xi_l}_{\xi_{l-1}}\psi(\sigma_1,\sigma_2)d\sigma_1 d\sigma_2+
$$
$$
+\sum\limits^{i-1}_{k=i-[(L-1)/2]}\sum\limits_{l=j+1}^{j+[(L-1)/2]}
h(\xi_i,\xi_j, \xi_{k-1},\xi_{l+1})\int\limits^{\xi_k}_{\xi_{k-1}}
\int\limits_{\xi_l}^{\xi_{l+1}}\psi(\sigma_1,\sigma_2)d\sigma_1 d\sigma_2+
$$
$$
\left.+\sum\limits^{i-1}_{k=i-[(L-1)/2]}\sum\limits^{j-1}_{l=j-[(L-1)/2]}
h(\xi_i,\xi_j, \xi_{k-1},\xi_{l-1})\int\limits^{\xi_k}_{\xi_{k-1}}
\int\limits^{\xi_l}_{\xi_{l-1}}\psi(\sigma_1,\sigma_2)d\sigma_1
d\sigma_2\right]+o\left(\frac{1}{n^r}\right) \ge
$$
$$
\ge o\left(\frac{1}{n^r}\right) +
$$
$$
+\frac{1}{L^2}\left[\sum\limits_{k=i+1}^{i+[(L-1)/2]}
\sum\limits_{l=j+1}^{j+[(L-1)/2]}\int\limits_{\xi_k}^{\xi_{k+1}}
\int\limits_{\xi_l}^{\xi_{l+1}}\psi(\sigma_1,\sigma_2)d\sigma_1
d\sigma_2\right.
\sum\limits^{L-1}_{i=0}\sum\limits^{L-1}_{j=0}h(\xi_i,\xi_j, 
\xi_{k+1},\xi_{l+1})+
$$
$$
+\sum\limits_{k=i+1}^{i+[(L-1)/2]}\sum\limits^{l=j-1}_{j-[(L-1)/2]}
\int\limits_{\xi_k}^{\xi_{k+1}}\int\limits^{\xi_l}_{\xi_{l-1}}
\psi(\sigma_1,\sigma_2)d\sigma_1 d\sigma_2\sum\limits^{L-1}_{i=0}
\sum\limits^{L-1}_{j=0}
h(\xi_i,\xi_j, \xi_{k+1},\xi_{l-1})+
$$
$$
+\sum\limits^{i-1}_{k=i-[(L-1)/2]}\sum\limits_{l=j+1}^{j+[(L-1)/2]}
\int\limits^{\xi_{k}}_{\xi_{k-1}}\int\limits_{\xi_l}^{\xi_{l+1}}
\psi(\sigma_1,\sigma_2)d\sigma_1 d\sigma_2
\sum\limits^{L-1}_{i=0}\sum\limits^{L-1}_{j=0}
h(\xi_i,\xi_j, \xi_{k-1},\xi_{l+1})+
$$
$$
 \left.+\sum\limits^{i-1}_{k=i-[(L-1)/2]}\sum\limits^{j-1}_{l=j-[(L-1)/2]}
\int\limits^{\xi_{k}}_{\xi_{k-1}}\int\limits^{\xi_l}_{\xi_{l-1}}
\psi(\sigma_1,\sigma_2)d\sigma_1 d\sigma_2
\sum\limits^{L-1}_{i=0}\sum\limits^{L-1}_{j=0}
h(\xi_i,\xi_j, \xi_{k-1},\xi_{l-1}) \right]=
$$
$$
=o(\frac{1}{n^r})+ 
\frac{1}{4\pi^2}\int\limits_0^{2\pi}\int\limits_0^{2\pi}
\psi(\sigma_1,\sigma_2)d\sigma_1 d\sigma_2\times
$$
$$
\times\left(\int\limits_0^{2\pi}\int\limits_0^{2\pi}
\frac{d\sigma_1 d\sigma_2}{(sin^2(\sigma_1/2)+sin^2(\sigma_2/2))^\lambda}
+O\left(
\left(\frac{logn}{n}\right)^{2-2\lambda}\right)\right),   \eqno (5.4)
$$
where the following relation was used:
$$
\frac{4\pi^2}{L^2}\sum\limits_{i=0}^{L-1}\sum\limits_{j=0}^{L-1}
 h(\xi_i,\xi_j,\xi_{k-1},\xi_{l-1})= O\left(\frac{\log n}{n}\right)+
 \int\limits_0^{2\pi}\int\limits_0^{2\pi}\frac{d\sigma_1d\sigma_2}
 {\left[sin^2(\sigma_1/2)+sin^2(\sigma_2/2)\right]^\lambda}.
$$

Without loss of generality one may assume
$k=1, l=1$ in the previous equation. Let us estimate
$$
U_0=\left|\frac{4\pi^2}{L^2}\sum\limits_{i=0}^{L-1}\sum\limits_{j=0}^{L-1}
 h(\xi_i,\xi_j,\ 0,0)-
 \int\limits_0^{2\pi}\int\limits_0^{2\pi}\frac{d\sigma_1d\sigma_2}
 {\left(sin^2(\sigma_1/2)+sin^2(\sigma_2/2)\right)^\lambda}\right|\le
$$
$$
\le\left|\sum\limits_{i=0}^{L-1}\sum\limits_{j=0}^{L-1}{}'
\int\limits_{\xi_i}^{\xi_{i+1}}
\int\limits_{\xi_j}^{\xi_{j+1}}\left[\frac{1}
 {\left(sin^2((\xi_i)/2)+sin^2((\xi_j)/2)\right)^\lambda}-\right.\right.
$$
$$
\left.\left.  -\frac{1}{\left(sin^2(\sigma_1/2)+sin^2(\sigma_2/2)\right)^\lambda}
\right] d\sigma_1d\sigma_2\right|+
$$
$$
+\left|\int\limits_{0}^{\xi_{1}}
\int\limits_{0}^{\xi_1}
\frac{1}{\left(sin^2(\sigma_1/2)+sin^2(\sigma_2/2)\right)^\lambda}
d\sigma_1d\sigma_2\right| = u_1+u_2,
$$
where $\sum\sum{}'$ means summation over  $(i,j)\neq (0,0).$

Let us estimate $u_1$ and $ u_2$.
One has:
$$
u_1\le\left|\sum\limits_{i=0}^{L-1}\sum\limits_{j=0}^{L-1}{}'
\int\limits_{\xi_i}^{\xi_{i+1}}
\int\limits_{\xi_j}^{\xi_{j+1}}\left[\frac{1}
 {\left(sin^2((\sigma_1)/2)+sin^2((\sigma_2)/2)\right)^\lambda}- \right.
 \right.
$$
$$
\left. \left. -\frac{1}{\left(sin^2(\xi_i/2)+sin^2(\sigma_2/2)\right)^\lambda}
\right] d\sigma_1d\sigma_2\right|+
$$
$$
+\left|\sum\limits_{i=0}^{L-1}\sum\limits_{j=0}^{L-1}{}'
\int\limits_{\xi_i}^{\xi_{i+1}}
\int\limits_{\xi_j}^{\xi_{j+1}}\left[\frac{1}
 {\left(sin^2((\xi_i)/2)+sin^2((\sigma_2)/2)\right)^\lambda}- \right.\right.
$$
$$
\left.\left. -\frac{1}{\left(sin^2(\xi_i/2)+sin^2(\xi_j/2)\right)^\lambda}
\right] d\sigma_1d\sigma_2\right|=
$$
$$
=u_{11}+u_{12}.
$$

The expressions $u_{11}$ and  $u_{12}$ can be estimated similarly. Let us
estimate $u_{11}$:
$$
u_{11}\le\frac c{L^4}\sum\limits_{i=0}^{L}\sum\limits_{j=0}^{L}{}'
\frac{1}{\left(sin^2((\xi_i)/2)+sin^2((\xi_j)/2)\right)^{1+\lambda}}\le
$$
$$
\le\frac c{L^{2-2\lambda}}\sum\limits_{i=0}^{L}\sum\limits_{j=0}^{L}{}'
\frac{1}{(i^2+j^2)^{1+\lambda}}\le c \frac1{L^{2-2\lambda}},
$$
where $c>0$ stands for various estimation constants.
Hence
$$
u_1\le \frac c{L^{2-2\lambda}}.
$$

Let us estimate $u_2:$
$$
u_2=\left|\int\limits_{0}^{\xi_{1}}
\int\limits_{0}^{\xi_1}
\frac{1}{\left(sin^2(\sigma_1/2)+sin^2(\sigma_2/2)\right)^\lambda}
d\sigma_d\sigma_2\right| \le
$$
$$
\le c \int\limits_{0}^{\xi_{1}}
\int\limits_{0}^{\xi_1}
\frac{1}{\left(\sigma_1^2+\sigma_2^2\right)^{\lambda}}
d\sigma_1d\sigma_2.
$$

Using polar coordinates, one gets:
$$
u_2 \le c\int\limits_{0}^{1/L}
\int\limits_{0}^{2\pi}
\frac{1}{\rho^{2\lambda-1}}d\rho d\phi \le \frac c{L^{2-2\lambda}}.
$$
Thus:
$$
U_0\le  \frac c{L^{2-2\lambda}}.
$$

From Lemmas 4.4, 4.1, and Theorem 4.1
it follows that
$$
\int\limits_0^{2\pi}\psi_1(\sigma_1)d\sigma_1 \ge
\frac{(1+o(1))(2\pi)^{r+1/q}R_{rq}(1)}{2^rr!(rq+1)^{1/q}(n-1+[R_{rq}(1)]^{1/r})^r},
\eqno (5.5)
$$
where
$R_{rq}(t)$ is a polynomial of degree $r$, least deviating from zero
in $L_q([-1,1]).$

Theorem 5.4 follows from inequalities  (5.4) and (5.5).

$\blacksquare$

{\bf 5.2. Optimal cubature formulas for calculating integrals (1.1).}

{\bf H\"older class of functions.}

Let $x_k:=2k\pi/n,$ $k=0,1,\ldots,n,$
$\Delta_{kl}=[x_k,x_{k+1},x_l,x_{l+1}],$ $k,l=0,1,\ldots,n-1,$
$x'_k=(x_{k+1}+x_k)/2,$ $k=0,1,\ldots,n-1,$ and
  $(s_1,s_2) \in \Delta_{ij},$ $i,j=0,1,\ldots,n-1.$

Calculate the integral  $Kf$ by the formula:
$$
Kf=\sum\limits^{n-1}_{k=0}\sum\limits^{n-1}_{l=0} f(x'_k,x'_l)
\int\int\limits_{\Delta_{kl}}
\frac{d\sigma_1 d\sigma_2}{\left(sin^2\left(\frac
{\sigma-x'_i}{2}\right)+
sin^2\left(\frac{\sigma-x'_j}{2}\right)\right)^\lambda}+
R_{nn}. \eqno (5.6)
$$

{\bf Theorem  5.5.} {\it Let $\Psi=H_{\alpha\alpha}(D), 0<\alpha < 1.$
Among all cubature formulas (5.1) with $\rho_1=\rho_2=0,$
formula (5.6), which has the error
$$
R_{nn}[\Psi]=\frac{(2+o(1))\gamma}{1+\alpha}\left(\frac{\pi}{n}\right)^\alpha,
$$
is asymptotically optimal. Here $\gamma$ is defined in (5.1').
}

{\bf Proof.}
Using the periodicity of the integrand, we estimate the error of
 cubature formula (5.6) as follows:
$$
|R_{nn}| \le \left|\sum\limits^{n-1}_{k=0}\sum\limits^{n-1}_{l=0}
\int\int\limits_{\Delta_{kl}}\left[\frac{f(\sigma_1,\sigma_2)-f(x'_i,x'_j)}
{\left(sin^2\frac{\sigma_1-s_1}{2}+
sin^2\frac{\sigma_2-s_2}{2}\right)^\lambda}-\right.\right.
$$
$$
-\left.\left.\frac{f(x'_k,x'_l)-f(x'_i,x'_j)}
{(\left(sin^2\frac{\sigma_1-x'_i}{2}+
sin^2\frac{\sigma_2-x'_j}{2}\right)^\lambda}\right]
d\sigma_1d\sigma_2\right|\le
$$
$$
\le \left|\sum\limits^{n-1}_{k=0}\sum\limits^{n-1}_{l=0}
\int\int\limits_{\Delta_{kl}}\frac{f(\sigma_1,\sigma_2)-f(x'_k,x'_l)}
{\left(sin^2\frac{\sigma_1-s_1}{2}+
sin^2\frac{\sigma_2-s_2}{2}\right)^\lambda} d\sigma_1d\sigma_2\right|+
$$
$$
+\left|\sum\limits^{n-1}_{k=0}\sum\limits^{n-1}_{l=0}
\int\int\limits_{\Delta_{kl}}(f(x'_k,x'_l)-f(x'_i,x'_j))\times  \right.
$$
$$
\left. \times\left[\frac1{\left(sin^2\frac{\sigma_1-s_1}{2}+
sin^2\frac{\sigma_2-s_2}{2}\right)^\lambda}-
\frac1{\left(sin^2\frac{\sigma_1-x'_i}{2}+
sin^2\frac{\sigma_2-x'_j}{2}\right)^\lambda}\right]
d\sigma_1d\sigma_2\right|=
$$
$$
=r_1+r_2.
$$

Let us estimate each of the sums $r_1$ and  $r_2$ separately. One has:
$$
r_1 \le \left|\sum\limits^{i+M}_{k=i-M}\sum\limits^{j+M}_{l=j-M}
\int\int\limits_{\Delta_{kl}}\left[\frac{f(\sigma_1,\sigma_2)-f(x'_k,x'_l)}
{\left(sin^2\frac{\sigma_1-s_1}{2}+
sin^2\frac{\sigma_2-s_2}{2}\right)^\lambda}d\sigma_1d\sigma_2\right|+ \right.
$$
$$
+ \left|\sum\limits^{n-1}_{k=0}\sum\limits^{n-1}_{l=0}{}'
\int\int\limits_{\Delta_{kl}}\left[\frac{f(\sigma_1,\sigma_2)-f(x'_k,x'_l)}
{\left(sin^2\frac{\sigma_1-s_1}{2}+
sin^2\frac{\sigma_2-s_2}{2}\right)^{\lambda}}d\sigma_1d\sigma_2\right|=\right.
$$
$$
=r_{11}+r_{12},
$$
where  $\sum\sum{}'$ means summation over $(k,l)$ such that \\
$\Delta_{kl}\notin \Delta^*,\  \Delta^*=[x_{i-M},x_{i+M+1};x_{j-M},x_{j+M+1}]
, M=[ln n].$

Furthermore
$$
r_{11} \le \frac c{n^\alpha}
\int\int\limits_{\Delta^*}\frac{d\sigma_1d\sigma_2}
{\left(sin^2\frac{\sigma_1-s_1}{2}+
sin^2\frac{\sigma_2-s_2}{2}\right)^\lambda}\le
$$
$$
 \le \frac c{n^\alpha}
\int\limits_0^{2\pi M/n}\int\limits_{0}^{2\pi}
\frac{d\rho d\phi}
{\rho^{2\lambda -1}} \le  \frac {c \log n}{n^{\alpha+2-2\lambda}}=
o\left(\frac 1{n^\alpha}\right).
$$

Estimating $r_{12},$ one can assume without loss of generality
$(i,j)=(0,0),$ and get:
$$
r_{12}\le 4\int\limits_0^{\pi/n}\int\limits_0^{\pi/n}
(\omega_1(\sigma_1)+\omega_2(\sigma_2))d\sigma_1d\sigma_2
\sum_{k=0}^{n-1}\sum_{l=0}^{n-1}h_{kl}(s_1,s_2,\sigma_1,\sigma_2)\le
$$
$$
\le 4\int\limits_0^{\pi/n}\int\limits_0^{\pi/n}
(\sigma_1^\alpha+\sigma_2^\alpha)d\sigma_1d\sigma_2
\sum_{k=0}^{n-1}\sum_{l=0}^{n-1}h_{kl}(s_1,s_2,\sigma_1,\sigma_2)\le
$$
$$
\le \frac 8{1+\alpha}\left(\frac \pi{n}\right)^{2+\alpha}
\sum_{k=0}^{n-1}\sum_{l=0}^{n-1}h_{kl}(s_1,s_2,\sigma_1,\sigma_2)\le
$$
$$
\le \frac{1+o(1)}{1+\alpha}2\left(\frac\pi{n}\right)^{\alpha}
\int\limits_0^{2\pi}\int\limits_0^{2\pi}
\frac{d\sigma_1d\sigma_2}
{\left(sin^2\frac{\sigma_1}{2}+
sin^2\frac{\sigma_2}{2}\right)^\lambda}.
$$

Here
$$h_{kl}(s_1,s_2;\sigma_1,\sigma_2)
=\sup_{(\sigma_1,\sigma_2)\in \Delta_{kl}}h(s_1,s_2;\sigma_1,\sigma_2).
$$

Combining the estimates of $r_{11}$ and  $r_{12}$, one gets:
$$
r_1\le \frac{1+o(1)}{1+\alpha}2\left(\frac\pi{n}\right)^{\alpha}\gamma
$$

Let us estimate $r_2.$ To this end we estimate the difference
$$
r_2(k,l)=\int\int\limits_{\Delta_{kl}}\left|
f(x'_k,x'_l)-f(x'_i,x'_j)
\right|\times
$$
$$
\times\left|\left[\frac 1{\left(sin^2\frac{\sigma_1-s_1}{2}+
sin^2\frac{\sigma_2-s_2}{2}\right)^\lambda}-
\frac{1}
{\left(sin^2\frac{\sigma_1-x'_i}{2}+
sin^2\frac{\sigma_2-x'_j}{2}\right)^\lambda}\right]
d\sigma_1d\sigma_2\right| .
$$

First, we estimate
$$
r_2(i,j)\leq \frac c{n^\alpha}\int\int\limits_{\Delta_{ij}}\left|
\frac1{\left(sin^2\frac{\sigma_1-s_1}{2}+
sin^2\frac{\sigma_2-s_2}{2}\right)^\lambda}-
\frac{1}
{\left(sin^2\frac{\sigma_1-x'_i}{2}+
sin^2\frac{\sigma_2-x'_j}{2}\right)^\lambda}\right|
d\sigma_1d\sigma_2\le
$$
$$
\le \frac c{n^{2+\alpha-2\lambda}}.
$$

The value  $r_2(k,l)$ is estimated similarly for $|k-i|\le 3$ and
$|l-j|\le 3.$

Let us estimate  $r_2(k,l)$ for other values of $k$ and
$l.$

One has:
$$
r_2(k,l)=\int\int\limits_{\Delta_{kl}}
\left|f(x'_k,x'_l)-f(x'_i,x'_j)\right|\times
$$
$$
\times\left|\frac1{\left(sin^2\frac{\sigma_1-s_1}{2}+
sin^2\frac{\sigma_2-s_2}{2}\right)^\lambda}-
\frac{1}
{\left(sin^2\frac{\sigma_1-x'_i}{2}+
sin^2\frac{\sigma_2-x'_j}{2}\right)^\lambda}\right|
d\sigma_1d\sigma_2\le
$$
$$
\le \frac{c}{n}\int\int\limits_{\Delta_{kl}}
\left[|x'_k-x'_i|^\alpha+|x'_l-x'_j|^\alpha\right]
\left[\left(\frac{|k-i|}n\right)+\left(\frac{|l-j|}n\right)\right]\times
$$
$$
\times\left|\frac1{\left(sin^2\frac{\sigma_1-x'_i+\theta_1(s_1-x'_i)}{2}+
sin^2\frac{\sigma_2-s_2}{2}\right)^{1+\lambda}}+ \right.
$$
$$
\left. +\frac{1}
{\left(sin^2\frac{\sigma_1-x'_i+\theta_1(s_1-x'_i)}{2}+
sin^2\frac{\sigma_2-x'_j+\theta_2(s_2-x'_j)}{2}\right)^{1+\lambda}}\right|
\le
$$
$$
\le \frac c{n^3}\left(\left|\frac{|k-i|}n\right|^\alpha+
\left|\frac{|l-j|}n\right|^\alpha\right)
\left(\left|\frac{|k-i|}n\right|+
\left|\frac{|l-j|}n\right|\right)
\left(\frac{n^2}{|k-i|^2+|l-j|^2}\right)^{1+\lambda}\le
$$
$$
\le \frac c{n^{\alpha+2-2\lambda}}\frac{(|k-i|+|l-j|)^{1+\alpha}}
{\left(|k-i|^2+|l-j|^2\right)^{1+\lambda}}\le
$$
$$
\le \frac c{n^{\alpha+2-2\lambda}}
\frac{(|k-i|^2+|l-j|^2)^{(1+\alpha)/2}}
{\left(|k-i|^2+|l-j|^2\right)^{1+\lambda}}\le
$$
$$
\le \frac c{n^{\alpha+2-2\lambda}}
\frac{1}
{\left(|k-i|^2+|l-j|^2\right)^{1/2-\alpha/2+\lambda}}.
$$

To estimate $r_2$, one sums up the last expression
over  $k$ and  $l.$ Without loss of generality assume
 $(i,j)=(0,0).$ Then
$$
r_2 \le \frac {c}{n^{\alpha+2-2\lambda}}\left( 16+4\sum_{k=0}^{[n/2]+1}
\sum_{l=0}^{[n/2]+1}{}'\frac 1{(k^2+l^2)^{\lambda+1/2-\alpha/2}}\right),
$$
where $\sum\sum'$ means summation over  $k$ and $l$ such that
$k>3$ or $l>3.$

One has:
$$
\sum_{k=0}^{[n/2]+1}
\sum_{l=0}^{[n/2]+1}{}'\frac 1{(k^2+l^2)^{\lambda+1/2-\alpha/2}}\le
$$
$$
\le A\left[\sum_{k=3}^{[n/l]+1}
\frac1{k^{2\lambda+1-\alpha}}+\sum_{k=3}^{[n/2]+1}
\sum_{l=3}^{[n/2]+1}\frac 1{(k^2+l^2)^{\lambda+1/2-\alpha/2}}\right]\le
$$
$$
\le A
\left \{
\begin{array}{ccc}
1, \quad if \quad 2\lambda - \alpha >1;\\
\log n, \quad if \quad 2\lambda - \alpha =1;\\
n^{1-2\lambda+\alpha},  \quad if \quad 2\lambda - \alpha <1.\\
\end{array}
\right.
$$

Hence
$$
{\bf r_2 \le }
 A
\left \{
\begin{array}{ccc}
n^{-(\alpha+2-2\lambda)},  \quad if \quad 2\lambda - \alpha >1;\\
n^{-1}\log n,  \quad if \quad  2\lambda - \alpha =1;\\
n^{-1},  \quad if \quad  2\lambda - \alpha <1.\\
\end{array}
\right.
$$

Thus, if $\alpha<1,$ then
$$
r_2\le o(n^{-\alpha}).
$$

Combining the estimates of $r_1$ and $r_2,$
one gets:
$$
R_{nn}[\Psi]\le \gamma
\frac{(2+o(1))}{1+\alpha}\left(\frac{\pi}{n}\right)^\alpha.
$$
Theorem 5.5 follows from the comparison of this inequality
with the lower bound of the value
$\zeta_{nn}[H_{\alpha,\alpha}(D)],$ mentioned in the  Corollary
to Theorem 5.1.
$\blacksquare$

{\bf Remark.} {\it If $\alpha=1,$ the cubature formula (5.6)
is optimal with respect to order.}

The proof of Theorem 5.5 yields also the following result:

{\bf Theorem 5.5'.}
{\it Let $\Psi=H_{\alpha\alpha}(D), 0<\alpha \le 1.$
Among all possible cubature formulas (5.1) with $\rho_1=\rho_2=0,$
formula
$$
Kf=\sum\limits^{n-1}_{k=0}\sum\limits^{n-1}_{l=0} f(x'_k,x'_l)
\int\int\limits_{\Delta_{kl}}
\frac{d\sigma_1 d\sigma_2}{\left(sin^2\left(\frac
{\sigma-s_1}{2}\right)+
sin^2\left(\frac{\sigma-s_2}{2}\right)\right)^\lambda}+
R_{nn},
$$
which has the error
$$
R_{nn}[\Psi]=\frac{(2+o(1))\gamma}{1+\alpha}\left(\frac{\pi}{n}\right)^\alpha,
$$
is asymptotically optimal.
}

To apply formula (5.6), one has to calculate
the integrals
$$
I_{kl}=\int\int\limits_{\Delta_{kl}}\frac{d\sigma_1d\sigma_2}
{\left(sin^2\frac{\sigma_1-x_i'}{2}+sin^2\frac{\sigma_2-x_j'}{2}
\right)^\lambda}
\eqno (5.7)
$$
for $k,l=0,1,\ldots,n-1.$ Exact values of these integrals
for arbitrary values  $\lambda$ are apparently unknown. Therefore
the procedure of numerical calculation of integrals (5.7)
should be given for practical application of formula (5.6).

Let $k=i$ and $l=j.$ Then the integral  $I_{ij}$ is replaced by the
integral
$$
p_{ij}{}^*=
\int\limits_{-\frac{\pi}{n}}^{\frac{\pi}{n}}\int\limits_{-\frac{\pi}{n}}^{\frac{\pi}{n}}
\frac{d\sigma_1d\sigma_2}
{\left(sin^2\frac{\sigma_1}{2}+sin^2\frac{\sigma_2}{2}\right)^{\lambda}+h},
\quad h>0,
$$
which can be calculated by cubature formulas ( in particular, Gauss
quadrature rule ) with arbitrary degree
of accuracy because the function

$
\frac{1}
{\left(sin^2\frac{\sigma_1}{2}+sin^2\frac{\sigma_2}{2}\right)^{\lambda}+h},
$

has derivatives up to arbitrary order.
The choice of parameter  $h$ is discussed in
Section 8.

Let $k=i, l \ne j,$ and
$$
I_{il}= \frac{4\pi^2}{n^2}\left(sin^2\frac{x'_l-x'_j}{2}\right)^{-\lambda}=
p^*_{il}.
$$

Let $k \ne i,$  $l=j,$ and
$$
I_{kj} =\frac{4\pi^2}{n^2}\left(sin^2\frac{x'_k-x'_i}{2}\right)^{-\lambda}=
p^*_{kj}.
$$

Let  $k \ne i,$ $l \ne j$, and
$$
I_{kl} = \frac{4\pi^2}{n^2}\left(sin^2\frac{x'_k-x'_i}{2}+sin^2\frac{x'_l-x'_j}{2} \right)^{-\lambda}=
p^*_{kl}.
$$

The integral  $Kf$ is calculated by the formula
$$
Kf=\sum\limits_{k=0}^{n-1}\sum\limits_{l=0}^{n-1}
p^*_{kl}f(x'_k,x'_l)+R_{nn}(f,p^*_{kl},x_k,y'_l). \eqno (5.8)
$$

Formula (5.8) is not optimal since it is not exact on
constant functions $f(x,y)=const.$ But one can estimate the error of this
formula:
$$
|R_{nn}(f,p^*_{kl},x'_k,y'_l))| \le M \sum\limits_{k=0}^{n-1}\sum\limits_{l=0}^{n-1}
|I_{kl}-p^*_{kl}|+R_{nn}(\Psi),
$$
where $M=\max|f(x,y)|.$

The values  $|I_{kl}-p^*_{kl}|$ are easily estimated, and one gets
the conclusion of Theorem 5.5'.

{\bf Classes of smooth functions}

{\bf Theorem 5.6.} {\it Assume $\varphi \in \tilde W^{r,r}(1).$
Let  $\Psi=\tilde W^{r,r}(1),$ and calculate the
integral $K\varphi$ by formula (5.1)
with $\rho_1=r-1$, $\rho_2=r-1$, and $n_1=n_2=n.$
Then the cubature formula
$$
K\varphi=\int\limits_{0}^{2\pi}\int\limits_{0}^{2\pi}
\frac{\varphi_{mn}(\sigma_1,
\sigma_2)d\sigma_1d\sigma_2}{(sin^2(\sigma_1-s_1)/2+ sin^2(\sigma_2-s_2)/2)^\lambda}+
R_{mn}(\varphi) \eqno (5.9)
$$
is asymptotically optimal. }

Before proving Theorem 5.6, let us describe the construction of the spline
$\varphi_{mn}.$
Let  $x_k=2k\pi/n, \ k=0,1,\dots,n.$
Divide the sides of the squares  $\Omega=[0,2\pi;0,2\pi]$ into  $n$
equal parts. Denote by
$\Delta_{kl}$ the rectangle $\Delta_{kl}=[2k\pi/n,2(k+1)\pi/n;
2l\pi/n, 2(l+1)\pi/n], k,l=0,1,\dots,n-1.$ Let
$(s_1,s_2)\in\Delta_{ij}.$ First we approximate
$\varphi(\sigma_1,\sigma_2)$ as a function of $\sigma_2,$ and
construct a spline  $\varphi_n(\sigma_1,\sigma_2)$ by the following rule.
Let  $\sigma_1$ be an arbitrary fixed number, $0\le \sigma_1\le 2\pi.$
On the segments  $[2k\pi/n, 2(k+1)\pi/n]$
for $k\ne j-2,\dots,j+1,$ one has:
$$\varphi_n(\sigma_1,\sigma_2)=\sum_{l=0}^{r-1} \left[\frac
{\varphi^{(0,l)}(\sigma_1,2k\pi/n)}{l!}(\sigma_2-2k\pi/n)^l+
B_l\delta^{(l)}(\sigma_1,(k+1)/n)\right],$$
where
$$\delta(\sigma_1,\sigma_2):=\varphi(\sigma_1,\sigma_2)-
\sum_{l=0}^{r-1}\frac {\varphi^{(0,l)}(\sigma_1,2k\pi/n)}{l!}
(\sigma_2-2k\pi/n)^l.
$$

The coefficients  $B_l$ are defined by the equation
$$
\left(2(k+1)\pi/n-\sigma_2\right)^r-\sum_{l=0}^{r-1}\frac {B_lr!}
{(r-l-1)!}\frac {2\pi}n\left(2\pi(k+1)/n-\sigma_2\right)^{r-l-1}=
$$
$$
=(-1)^rR_{r1}
\left(2\pi(2k+1)/2n; \pi/n; \sigma_2\right),$$
where $R_{r1}(a,h,x) $ is a polynomial of degree $r$, least deviating
from zero in the norm of the space $L$ on the segment $[a-h, a+h].$
On the segment $\left[2\pi(j-2)/n,2\pi(j+2)/n\right]$ the function $\varphi_n
(\sigma_1,\sigma_2)$ is defined by the partial sum of the Taylor series:
$$\varphi_n(\sigma_1,\sigma_2)=\varphi(\sigma_1,2\pi j/n)+
{\displaystyle\frac
{\varphi^{(0,1)}(\sigma_1,2\pi j/n)}{1!}}(\sigma_2-j/n)+
\cdots +$$
$$+{\displaystyle
 \frac {\varphi^{(0,r-1)}(\sigma_1,2\pi j/n)}{(r-1)!}}
(\sigma_2-2\pi j/n)^{r-1}.$$
We define the function $\varphi_{nn}(\sigma_1,\sigma_2)$ by analogy with
the function $\varphi_n(\sigma_1,\sigma_2).$

{\bf Proof of Theorem 5.6.} Let $(s_1,s_2)\in \Delta_{ij}.$
The error of formula (5.9) we estimate by the inequality
$$|R_{nn}|\leqslant\sum_{k=0}^{n-1}\sideset{}{'}\sum_{l=0}^{n-1}\left|
\iint\limits_{\Delta_{kl}}\frac{\varphi(\sigma_1,\sigma_2)-
\varphi_{nn}(\sigma_1,\sigma_2)}
{\Bigl(\sin^2\frac{\sigma_1-s_1}{2}+
\sin^2\frac{\sigma_2-s_2}{2}\Bigr)^{\lambda}}d\sigma_1d\sigma_2\right|+
$$
%------------------
$$
  +\sum_{k=0}^{n-1}\sideset{}{''}\sum_{l=0}^{n-1}\left|
  \int\limits_{\Delta_{kl}}\frac{\varphi(\sigma_1,\sigma_2)-
  \varphi_{nn}(\sigma_1,\sigma_2)}
  {\Bigl(\sin^2\frac{\sigma_1-s_1}{2}+
  \sin^2\frac{\sigma_2-s_2}{2}\Bigr)^{\lambda}}d\sigma_1d\sigma_2\right|=r_1+r_2,
\eqno (5.10)
$$
where $\sideset{}{'}\sum\limits_{k,l}$ means summation over $(k,l)$
such that $i-1\leqslant k \leqslant i+1,$ $0\leqslant l \leqslant
n-1$ or $0\leqslant k \leqslant n-1,$ $j-1\leqslant l \leqslant
j+1,$ and $\sideset{}{''}\sum\limits_{k,l}$ means summation over
the other values of $(k,l).$

Let us estimate each of the sums $r_1$ and $r_2$ separately. In
addition without loss of generality assume
that $ \iint\limits_{\Delta_{kl}}(\varphi(\sigma_1,\sigma_2)-
  \varphi_{nn}(\sigma_1,\sigma_2))d\sigma_1d\sigma_2 \geqslant 0.$
Then
$$r_1\leqslant\sum_{k=0}^{n-1}\sideset{}{'}\sum_{l=0}^{n-1}
|\varphi(\sigma_1,\sigma_2)- \varphi_{nn}(\sigma_1,\sigma_2)|
\iint\limits_{\Delta_{kl}}\frac{d\sigma_1d\sigma_2}
{\Bigl(\sin^2\frac{\sigma_1-s_1}{2}+
\sin^2\frac{\sigma_2-s_2}{2}\Bigr)^{\lambda}}\leqslant $$
%-------------------
$$
  \leqslant A
  \begin{cases}
    n^{-(r+1)} &, \lambda \le 1/2 \\
    n^{-(r+2-2\lambda)} &, \lambda>1/2;
  \end{cases}
\eqno (5.11)
$$
%-------------------
$$r_2\leqslant 4\sum_{k=i+2}^{i+1+[(n-1)/2]}
\sum_{l=j+2}^{j+1+[(n-1)/2]}\frac{1}{\Bigl(\sin^2\frac{x_k-s_1}{2}+
\sin^2\frac{x_l-s_2}{2}\Bigr)^{\lambda}}
\iint\limits_{\Delta_{kl}}\psi(\sigma_1,\sigma_2)d\sigma_1d\sigma_2-
$$
%-------------------
$$-4\sum_{k=i+2}^{i+1+[(n-1)/2]}\sum_{l=j+2}^{j+1+[(n-1)/2]}
\Biggl[\frac{1}{\Bigl(\sin^2\frac{x_k-s_1}{2}+
\sin^2\frac{x_l-s_2}{2}\Bigr)^{\lambda}}-
\frac{1}{\Bigl(\sin^2\frac{x_{k+1}-s_1}{2}+
\sin^2\frac{x_{l+1}-s_2}{2}\Bigr)^{\lambda}}\Biggr]\cdot $$
%-------------------
$$
  \cdot\iint\limits_{\Delta_{kl}}
  \psi^-(\sigma_1,\sigma_2)d\sigma_1d\sigma_2
  =r_{21}+r_{22},
\eqno (5.12)
$$
where $\psi(\sigma_1,\sigma_2)=\varphi(\sigma_1,\sigma_2)-
\varphi_{nn}(\sigma_1,\sigma_2),$
%-------------------
$$ \psi^+(\sigma_1,\sigma_2)=
  \begin{cases}
    \psi(\sigma_1,\sigma_2) &, \text{ if } \psi(\sigma_1,\sigma_2)\geqslant0 \\
    0 &, \text{ if }\psi(\sigma_1,\sigma_2)<0;
  \end{cases}
$$
%--------------------
$$ \psi^-(\sigma_1,\sigma_2)=
  \begin{cases}
     0 &, \text{ if }\psi(\sigma_1,\sigma_2)\geqslant0;\\
    -\psi(\sigma_1,\sigma_2) &, \text{ if } \psi(\sigma_1,\sigma_2)<0.
  \end{cases}
$$
%--------------------

It is obvious
%--------------------
$$
 4\sum_{k=i+2}^{i+1+[(n-1)/2]}\sum_{l=j+2}^{j+1+[(n-1)/2]}
 \frac{1}{\Bigl(\sin^2\frac{x_k-s_1}{2}+
 \sin^2\frac{x_l-s_2}{2}\Bigr)^{\lambda}}\leqslant
 \frac{1+o(1)}{4\pi^2}\int\limits_{0}^{2\pi}\int\limits_{0}^{2\pi}
 \frac{d\sigma_1d\sigma_2}{\Bigl(\sin^2\frac{\sigma_1}{2}+
 \sin^2\frac{\sigma_2}{2}\Bigr)^{\lambda}}
\eqno (5.13)
$$
%-------------------

Let us estimate the integral
%-------------------
$$i=\iint\limits_{\Delta_{kl}}
\psi(\sigma_1,\sigma_2)d\sigma_1d\sigma_2\leqslant
\left|\iint\limits_{\Delta_{kl}}\Bigl(\varphi(\sigma_1,\sigma_2)-
\varphi_{n}(\sigma_1,\sigma_2)\Bigr)d\sigma_1d\sigma_2\right|+$$
%-------------------
$$
   +\left|\iint\limits_{\Delta_{kl}}\Bigl(\varphi_n(\sigma_1,\sigma_2)-
   \varphi_{nn}(\sigma_1,\sigma_2)\Bigr)d\sigma_1d\sigma_2\right|=i_1+i_2.
\eqno (5.14)
$$

Since the expressions $i_1$ and $i_2$ are estimated similarly, we
estimate only $i_1.$ One has:
%--------------------
$$ i_1\leqslant \frac{2\pi}{n}\max_{s_1}\left|\int\limits_{x_l}^{x_{l+1}}
\Bigl(\varphi(s_1,\sigma_2)-\varphi_{n}(s_1,\sigma_2)\Bigr)
d\sigma_2 \right|. $$
%--------------------
This integral is a continuous function of $s_1,$ which
attains its maximum at a point $s^*$, and
%--------------------
$$ i_1\leqslant \frac{2\pi}{n}\left|\int\limits_{x_l}^{x_{l+1}}
\Bigl(\varphi(s^*,\sigma_2)-\varphi_{n}(s^*,\sigma_2)\Bigr)
d\sigma_2 \right|\leqslant $$
%--------------------
$$\leqslant \frac{2\pi}{r!n}\int\limits_{x_l}^{x_{l+1}}
\bigl|\varphi^{(0,r)}(s^*,\sigma_2)\bigr|\Biggl|(x_{l+1}-\sigma_2)^r-
$$
%--------------------
$$-\sum_{j=0}^{r-1}\frac{B_{lj}(x_{l+1}-x_l)r!}
{(r-1-j)!}(x_{l+1}-\sigma_2)^{r-j-1}\Biggr| d\sigma_2\leqslant $$
%--------------------
$$\leqslant \frac{2\pi}{r!n}\int\limits_{x_l}^{x_{l+1}}
\Biggl|(x_{l+1}-\sigma_2)^r-\sum_{j=0}^{r-1}\frac{B_{lj}(x_{l+1}-x_l)r!}
{(r-1-j)!}(x_{l+1}-\sigma_2)^{r-j-1}\Biggr| d\sigma_2= $$
%--------------------
$$
   =\frac{2\pi}{r!n}\int\limits_{x_l}^{x_{l+1}}
   \bigl|R_{r1}(\sigma_2)\bigr| d\sigma_2\leqslant
   \frac{4}{(r+1)!}\left(\frac{\pi}{n}\right)^{r+2}R_{r1}(1).
\eqno (5.15)
$$

From inequalities (5.14) and (5.15) one gets:
$$i\leqslant
\frac{8}{(r+1)!}\left(\frac{\pi}{n}\right)^{r+2}R_{r1}(1)$$ and
%--------------------
$$
   r_{21}\leqslant
   \frac{2+o(1)}{(r+1)!}\left(\frac{\pi}{n}\right)^{r}R_{r1}(1)
   \int\limits_{0}^{2\pi}\int\limits_{0}^{2\pi}
   \frac{d\sigma_1d\sigma_2}{\Bigl(\sin^2\frac{\sigma_1}{2}+
   \sin^2\frac{\sigma_2}{2}\Bigr)^{\lambda}}.
\eqno (5.16)
$$

One has:
$$
   r_{22}=o(n^{-r}).
\eqno (5.17)
$$

Estimate (5.17) follows from the inequalities:
%--------------------
$$\left|\iint\limits_{\Delta_{kl}}
\psi^-(\sigma_1,\sigma_2)d\sigma_1d\sigma_2\right|\leqslant
\iint\limits_{\Delta_{kl}}
\bigl|\psi(\sigma_1,\sigma_2)\bigr|d\sigma_1d\sigma_2= O(n^{-r-2})$$
and
%-------------------
$$\sum_{k=i+2}^{i+1+[(n-1)/2]}\sum_{l=j+2}^{j+1+[(n-1)/2]}
\Biggl|\frac{1}{\Bigl(\sin^2\frac{x_k-s_1}{2}+
\sin^2\frac{x_l-s_2}{2}\Bigr)^{\lambda}}-
\frac{1}{\Bigl(\sin^2\frac{x_{k+1}-s_1}{2}+
\sin^2\frac{x_{l+1}-s_2}{2}\Bigr)^{\lambda}}\Biggr|\leqslant$$
%--------------------
$$\leqslant An^{2\lambda}\sum_{k}\sum_{l}\frac{(k-i)+(l-j)}
{\bigl((k-i)^2+(l-j)^2\bigr)^{\lambda+1}}\leqslant c
  \begin{cases}
    n &, \lambda<1/2 \\
    n\log n &, \lambda=1/2 \\
    n^{2\lambda} &, \lambda>1/2.
  \end{cases}
$$

The estimate
$$
R_{nn}(\Psi) \le
 (1+o(1))
\frac{2\pi^{r}R_{r1}(1)}{(r+1)!(n-1+[R_{r1}(1)]^{1/r})^r}
\gamma
$$
 follows from inequalities (5.10),
(5.11), (5.16), and (5.17).

Theorem 5.5 follows from the comparison of the values $\zeta_{nn}[\Psi]$
and $R_{nn}[\Psi].$
$\blacksquare$

Let us construct cubature formulas for calculating integrals
$Kf$ on  classes of functions $W^{rr}(1)$. These formulas will be
less accurate than the ones in Theorem 5.3, but they will be
optimal with respect to order,
and easier to apply.

First, we investigate the smooth  function
$$
\psi(t_1,t_2)=\int\limits_0^{2\pi}\int\limits_0^{2\pi}\frac
{f(\tau_1,\tau_2)d\tau_1 d\tau_2}{\left(\sin^2\frac{\tau_1-t_1}2+
\sin^2\frac{\tau_2-t_2}2\right)^\lambda}
$$
assuming $f(t_1,t_2)\in \tilde W^{r,r}.$
Changing the variables $\tau_1=\tau_1-t,\
\tau_2=\tau_2-t,$ in the last integral, one gets:
$$
\psi(t_1,t_2)=\int\limits_0^{2\pi}\int\limits_0^{2\pi}\frac
{f(\tau_1+t_1,\tau_2+t_2)d\tau_1 d\tau_2}{\left(\sin^2\frac{\tau_1}2+
\sin^2\frac{\tau_2}2\right)^\lambda}
$$

Thus, $\psi(t_1,t_2)\in W^{r,r}.$

{\bf Remark.} {\it It is known [9] that Kolmogorov and Babenko widths on the class
of functions $W^{r,r}(1)$ are equal to  $\delta_n(W^{r,r}(1))\asymp d_n
(W^{r,r}(1),C)\asymp \frac 1{n^{r/2}}.$
Hence the recovery of the function $\psi(t_1,t_2)$ using $n$ functionals
is not possible with accuracy greater than
$O\left(\frac 1{n^{r/2}}\right).$ More
precise conclusions are obtained in Theorems
5.3 and 5.4.}

Thus, for recovery of a function $\psi(t_1,t_2),\
(t_1,t_2)\in [0,2\pi]^2$ with the accuracy $O(n^{-r/2}),$ it is sufficient
to calculate the value of the function $\psi(t_1,t_2)$ at the nodes
 $(v_k,v_l),$ where $v_k=2k\pi/N,$ $k,l=0,1,\dots,N,$ and $N^2=n,$
and to use the local spline  $\psi_N(t_1,t_2)$ of degree $r$
with respect to each variable.

Let us describe the construction of such spline.

Assume for simplicity that $M:=N/r $ is an integer, and
 cover the domain  $[0,2\pi]^2$ with the
squares $\Delta_{kl}=
[w_k,w_l],$ $k,l=0,1,\dots,M-1,$ here $w_k=2k\pi/M,$ $k=0,\dots,M.$
Approximate the
function $\psi(t_1,t_2)$ in each
domain $\Delta_{kl}$ by the interpolation polynomial $\psi_N(t_1,t_2,
\Delta_{kl})$ constructed  on the nodes  $(x_i^k,x_j^l),\ i,j=0,1,
\dots,r,$ $x_i^k=w_k+\frac{2\pi}{Mr}i,$ $i=0,1,\dots,r.$

Denote the local spline, which is defined by the polynomials
$\psi_N(t_1,t_2,
\Delta_{kl})$, by  $\psi_N(t_1,t_2).$

If the values $\psi(v_k,v_l)$ are calculated by formula
(5.9) with the accuracy
$O(n^{-r/2}),$ then
$$
\|\psi(t_1,t_2)-\psi_N(t_1,t_2)\|_C\le O(n^{-r/2}).
$$
Therefore the spline  $\psi_N(t_1,t_2)$ is optimal with respect to
order, and a method for recovery of  the function
$\psi(t_1,t_2),$
which has the error  $O(n^{-r/2})$ (in the $\sup-$norm) is constructed.

{\bf 6. Optimal methods for calculating integrals of the form  $Tf.$}

{\bf Lower bounds for the functionals $\zeta_{mn}$ and $\zeta_N.$}

First we get a lower bound for the error of formula (2.1) with
 $\rho_1=\rho_2=0$ and $n_1=n_2=n$, on H\"older classes.

{\bf Theorem 6.1.} {\it  Let $\Psi=H_{\alpha\alpha}(D),$ and calculate the
integral  $Tf$
by formula (2.1) with $n_1=n_2=n$ and $\rho_1=\rho_2=0$.
Then the estimate:
$$
\zeta_{nn}[\Psi] \ge \frac{(1+o(1))}{2^{2\lambda}(1+\alpha)n^\alpha}
\int\limits_{-1}^1\int\limits_{-1}^1\frac{dt_1dt_2}{(\tau^2_1+t_1^2)^\lambda}
\eqno (6.1)
$$
holds.             }

{\bf Proof.}
Let  $n>0$ be an integer, $L=[n/\log n].$ Let
$v_k:=-1+2k/L,$
$k=0,1,\ldots,L.$ By $(\xi_k,\eta_l)$ we denote a set which is the  union
of nodes  $(x_i,y_j)$, $i,j=1,2,\ldots,n$ of formula
(2.1) and the nodes $(v_i,v_j)$, $i,j=1,2,\ldots,L.$
Let  $\Delta_{kl}=[v_k,v_{k+1};v_l,v_{l+1}], \ k,l=0,1,\dots,L-1.$
 Let  $0\leq \psi(t_1,t_2)\in H_{\alpha\alpha}(D),$
where $D=[-1,1]^2,$  vanishing at the nodes $(\xi_k,\eta_l)$,
$k,l=0,1,\ldots,N,$ where $N=n+L.$

Consider the integral
$$
(T\psi)(v_i,v_j)=\int\limits_{-1}^1\int\limits_{-1}^1
\frac{\psi(\tau_1,\tau_2)d\tau_1d\tau_2}{((\tau_1-v_i)^2+
(\tau_2-v_j)^2)^\lambda}=
$$
$$
=\left(\sum\limits_{k=i}^{L-1}\sum\limits_{l=j}^{L-1}+\sum\limits_{k=i}^{L-1}\sum\limits_{l=0}^{j-1}+
\sum\limits_{k=0}^{i-1}\sum\limits_{l=j}^{L-1}+\sum\limits_{k=0}^{i-1}\sum\limits_{l=0}^{j-1}\right) \times
$$
$$
\times\int\int\limits_{\Delta_{kl}}\frac{\psi(\tau_1,\tau_2)d\tau_1d\tau_2}
{((\tau_1-v_i)^2+(\tau_2-v_j)^2)^\lambda}\ge
$$
$$
\ge\sum\limits_{k=0}^{L-i-1}\sum\limits_{l=0}^{L-j-1}\left(\frac{L}{2}\right)^{2\lambda}
\frac{1}{((k+1)^2+(l+1)^2)^\lambda}\int\int\limits_{\Delta_{k+i,l+j}}
\psi(\tau_1,\tau_2)d\tau_1d\tau_2+
$$
$$
+\sum\limits_{k=0}^{L-i-1}\sum\limits_{l=0}^{j-1}\left(\frac{L}{2}\right)^{2\lambda}
\frac{1}{((k+1)^2+(l+1)^2)^\lambda}\int\int\limits_{\Delta_{k+i,j-l-1}}
\psi(\tau_1,\tau_2)d\tau_1d\tau_2+
$$
$$
+\sum\limits_{k=0}^{i-1}\sum\limits_{l=0}^{L-j-1}\left(\frac{L}{2}\right)^{2\lambda}
\frac{1}{((k+1)^2+(l+1)^2)^\lambda}\int\int\limits_{\Delta_{i-k-1,j+l}}
\psi(\tau_1,\tau_2)d\tau_1d\tau_2+
$$
$$
+\sum\limits_{k=0}^{i-1}\sum\limits_{l=0}^{j-1}\left(\frac{L}{2}\right)^{2\lambda}
\frac{1}{((k+1)^2+(l+1)^2)^\lambda}\int\int\limits_{\Delta_{i-k-1,j-l-1}}
\psi(\tau_1,\tau_2)d\tau_1d\tau_2=
$$
$$
=\sum\limits_{k=0}^{L-1}\sum\limits_{l=0}^{L-1}\left(\frac{L}{2}\right)^{2\lambda}
\frac{U(L-i-1-k)U(L-j-1-l)}{((k+1)^2+(l+1)^2)^\lambda}\int\int\limits_{\Delta_{k+i,l+j}}
\psi(\tau_1,\tau_2)d\tau_1d\tau_2+
$$
$$
+\sum\limits_{k=0}^{L-1}\sum\limits_{l=0}^{L-1}\left(\frac{L}{2}\right)^{2\lambda}
\frac{U(L-i-1-k)U(j-1-l)}{((k+1)^2+(l+1)^2)^\lambda}\int\int\limits_{\Delta_{k+i,j-l-1}}
\psi(\tau_1,\tau_2)d\tau_1d\tau_2+
$$
$$
+\sum\limits_{k=0}^{L-1}\sum\limits_{l=0}^{L-1}\left(\frac{L}{2}\right)^{2\lambda}
\frac{U(i-1-k)U(L-j-1-l)}{((k+1)^2+(l+1)^2)^\lambda}\int\int\limits_{\Delta_{i-k-1,j+l}}
\psi(\tau_1,\tau_2)d\tau_1d\tau_2+
$$
$$
+\sum\limits_{k=0}^{L-1}\sum\limits_{l=0}^{L-1}\left(\frac{L}{2}\right)^{2\lambda}
\frac{U(i-1-k)U(j-1-l)}{((k+1)^2+(l+1)^2)^\lambda}
\int\int\limits_{\Delta_{i-k-1,j-l-1}}
\psi(\tau_1,\tau_2)d\tau_1d\tau_2.
$$

Here $U(k)=1$ for $k \ge 0,$ and $U(k)=0$ for $k < 0.$

Averaging the above inequality over all  $i$ and $j,$
$i,j=0,1,\ldots,L-1,$ one gets:
$$
R_{nn}(\Psi,p_{kl};x_k,y_l) \ge \frac{1}{L^2}\sum\limits_{i=0}^{L-1}
\sum\limits_{j=0}^{L-1}
T(\psi)(\xi_i,\eta_j) \ge
$$
$$
\ge\frac{1}{L^{2-2\lambda}2^{2\lambda}}\sum\limits_{k=0}^{L-1}\sum\limits_{l=0}^{L-1}
\frac{1}{((k+1)^2+(l+1)^2)^\lambda}
\left[\sum\limits_{i=0}^{L-1}\sum\limits_{j=0}^{L-1}
U(L-i-1-k) \times\right.
$$
$$
\left.\times U(L-j-1-l)\int\int\limits_{\Delta_{k+i,l+j}}
\psi(\tau_1,\tau_2)d\tau_1d\tau_2+\right.
$$
$$
+\sum\limits_{i=0}^{L-1}\sum\limits_{j=0}^{L-1}
U(L-i-1-k)U(j-1-l)\int\int\limits_{\Delta_{k+i,j-l-1}}
\psi(\tau_1,\tau_2)d\tau_1d\tau_2+
$$
$$
+\sum\limits_{i=0}^{L-1} \sum\limits_{j=0}^{L-1}
U(i-1-k)U(L-j-1-l)\int\int\limits_{\Delta_{i-k-1,j+l}}
\psi(\tau_1,\tau_2)d\tau_1d\tau_2+
$$
$$
\left.+\sum\limits_{i=0}^{L-1}\sum\limits_{j=0}^{L-1}
U(i-1-k)U(j-1-l)\int\int\limits_{\Delta_{i-k-1,j-l-1}}
\psi(\tau_1,\tau_2)d\tau_1d\tau_2\right] \ge
$$
$$
\ge\frac{1}{L^{2-2\lambda}2^{2\lambda}}\sum\limits_{k=0}^{L-1}\sum\limits_{l=0}^{L-1}
\frac{1}{((k+1)^2+(l+1)^2)^\lambda}\left[\int\limits_{v_k}^1\int\limits_{v_l}^1
\psi(\tau_1,\tau_2)d\tau_1d\tau_2+\right.
$$
$$
\left.+\int\limits_{v_k}^1\int\limits_{-1}^{v_{L-l-2}}
\psi(\tau_1,\tau_2)d\tau_1d\tau_2+\int\limits^{v_{L-k-2}}_{-1}\int\limits_{v_l}^1
\psi(\tau_1,\tau_2)d\tau_1d\tau_2+\int\limits^{v_{L-k-2}}_{-1}\int\limits_{-1}^{v_{L-l-2}}
\psi(\tau_1,\tau_2)d\tau_1d\tau_2\right]\ge
$$
$$
\ge \frac{1}{L^{2-2\lambda}2^{2\lambda}}\sum\limits_{k=0}^{L-1}
\sum\limits_{l=0}^{L-1}
\frac{1}{((k+1)^2+(l+1)^2)^\lambda}\int\limits_{-1}^1\int\limits_{-1}^1
\psi(\tau_1,\tau_2)d\tau_1d\tau_2.
\eqno (6.2)
$$

From inequality  (6.2) it follows that
$$
\zeta_{nn}[H_{\alpha\alpha}(D)] \ge (1+o(1))\frac{1}{L^{2-2\lambda}2^{2\lambda}}
\sum\limits_{k=1}^{L-1}\sum\limits_{l=1}^{L-1}
\frac{1}{(k^2+l^2)^\lambda}
\int\limits_{-1}^1\int\limits_{-1}^1\psi(\tau_1,\tau_2)d\tau_1d\tau_2 =
$$
$$
=\frac{1+o(1)}{2^{2\lambda}4}
\int\limits_{-1}^1\int\limits_{-1}^1\frac{dt_1dt_2}{(t_1^2+t_2^2)^\lambda}
\int\limits_{-1}^1\int\limits_{-1}^1\psi(\tau_1,\tau_2)d\tau_1d\tau_2.
\eqno (6.3)
$$

From Theorem 4.2 and Lemma 4.4 it follows that
the inequality
$$
\int\limits_{-1}^1\int\limits_{-1}^1\psi(\tau_1,\tau_2)d\tau_1d\tau_2 \ge
\frac{4}{1+\alpha}\frac{1}{n^\alpha} \eqno (6.4)
$$
is valid for an arbitrary vector of the weights and the nodes $(X,Y,P)$ on
the class
 $H_{\alpha\alpha}(D).$

Theorem 6.1 follows from inequalities (6.3) and (6.4).
$\blacksquare$

{\bf Theorem 6.2.}
{\it Let  $\Psi= C_2^r(1),$ and calculate the integral $Tf$
by formula  (2.1) with $\rho_1=\rho_2=0.$ If $n_1=n_2=n,$ then
$$
\zeta_{nn}[\Psi]\ge(1+o(1))\frac{2K_r}{2^{2\lambda}(\pi n)^r}
\int\limits_{-1}^{1}\int\limits_{-1}^{1}
\frac{ds_1ds_2}{(s_1^2)+s_2^2))^\lambda},
$$
where $K_r$ is the Favard constant.}

{\bf Proof.} Let
$$
\psi(s_1,s_2)=\psi_1(s_1)+\psi_2(s_2),
$$
where $0\leq \psi_1(s)\in W^r(1),$
vanishes at the nodes $x_k,$ $k=1,2,\ldots,n,$ and  $0\leq \psi_2(s)\in
W^r(1)$ vanishes at the nodes $y_k,$ $k=1,2,\ldots,n.$

For arbitrary nodes $x_k,$ $k=1,2,\ldots,n,$ one has (see [11]):

$$
\int\limits_{-1}^{1}\psi_i(s)ds \ge \frac{2 K_r}{(\pi n)^r}, \  i=1,2.
$$
Thus the inequality
$$
\int\limits_{-1}^{1}\int\limits_{-1}^{1} \psi(s_1,s_2)ds_1ds_2 \ge
\frac{8K_r}{(\pi n)^r}
$$
holds for arbitrary nodes $(x_1,\ldots,x_n)$ and
$(y_1,\ldots,y_n).$

Theorem 6.2 follows from
this estimate and inequality  (6.3).
$\blacksquare$

{\bf Theorem 6.3.}
{\it Let $\Psi=W_p^{r,r}(1),$ $r=1,2,\ldots,$ $1\le p \le \infty,$
and calculate the integral $Tf$ by formula (2.1)
with
$\rho_1=\rho_2=r-1$ and $n_1=n_2=n.$
Then the estimate
$$
\zeta_{nn}[\Psi] \ge (1+o(1))
\frac{2^{1/q}R_{rq}(1)}{2^{2\lambda}r!(rq+1)^{1/q}(n-1+[R_{rq}(1)]^{1/r})^r}
\int\limits_{-1}^1\int\limits_{-1}^1
\frac{ds_1ds_2}{(s_1^2 +s_2^2)^\lambda},
\eqno (6.5)
$$
holds,
where $R_{rq}(t)$ is a polynomial of degree $r$, least deviating from zero
in $L_q([-1,1]).$}

{\bf Proof.} Let  $L=[n/\log n].$  Consider the nodes
$(v_k,v_l),$ $v_k=\frac{2 k}{L},$
$k,l=0,1,\ldots,L-1.$
By $(\xi_i,\eta_j),$ $i,j=0,1,\ldots,N-1,$ $N=n+L$
denote the union of the nodes  $(x_k,y_l)$ and $(\xi_i ,\xi_j).$
Let $\psi(s_1,s_2)=\psi_1(s_1)+\psi_2(s_2),$
where $0\leq \psi_1(s)\in W^r_p(1)$
vanishes with its derivatives up to order  $r-1$
at the nodes  $\xi_i,$ $i=0,1,\ldots,N-1,$ and
 $0\leq \psi_2(s)\in W^r_p(1)$
vanishes  with its derivatives up to order  $r-1$
at the nodes $\eta_j,$ $j=0,1,\ldots,N-1.$ Assume
that
$
\int\limits_{v_i}^{v_{i+1}}\psi_1(s)ds>0, \,\,\,
i=0,1,\ldots,N-1,
$ and
$
\int\limits_{v_j}^{v_{j+1}}\psi_2(s)ds>0, \,\,\,
j=0,1,\ldots,N-1.
$

Using the argument similar to the one in the proof of Theorem 6.1, one
gets:
$$
\zeta_{nn}(\Psi,p_{kl};v_k,v_l) \ge \frac{1}{L^2}\sum\limits_{i=0}^{L-1}
\sum\limits_{j=0}^{L-1}
T(\psi)(v_i,v_j) \ge
$$
$$
\ge \frac{1}{L^{2-2\lambda}2^{2\lambda}}\sum\limits_{k=0}^{L-1}
\sum\limits_{l=0}^{L-1}
\frac{1}{((k+1)^2+(l+1)^2)^\lambda}\int\limits_{-1}^1\int\limits_{-1}^1
\psi(\tau_1,\tau_2)d\tau_1d\tau_2=
$$
$$
=\frac{1+o(1)}{2^{2\lambda}4}
\int\limits_{-1}^1\int\limits_{-1}^1\frac{dt_1dt_2}{(t_1^2+t_2^2)^\lambda}
\int\limits_{-1}^1\int\limits_{-1}^1\psi(\tau_1,\tau_2)d\tau_1d\tau_2.
\eqno (6.6)
$$

From Theorem 4.1 and  lemma 4.4 it follows that
the inequality
$$
\int\limits_{-1}^1\int\limits_{-1}^1\psi(\tau_1,\tau_2)d\tau_1d\tau_2 \ge
(1+o(1))\frac{2^{2+1/q}R_{rq}(1)}{r!(rq+1)^{1/q}(n-1+[R_{rq}(1)]^{1/q})^r}
 \eqno (6.7)
$$
is valid for arbitrary weights and the nodes $(X,Y,P)$ on the class
 $H_{\alpha\alpha}(D).$

Theorem 6.3 follows from inequalities (6.6)- (6.7).
$\blacksquare$

{\bf Cubature formulas.}

Let us construct a cubature formula for calculating the integral $Tf$
on the H\"older class $H_{\alpha\alpha}(D)$. Let
$x_k:=-1+2k/n,$ $k=0,1,\ldots,n,$ $x'_k=(x_{k+1}+x_k)/2,$
$k=0,1,\ldots,n-1,$ and
$\Delta_{kl}=[x_k,x_{k+1};x_l,x_{l+1}],$ $k,l=0,1,\ldots,n-1.$

Calculate the integral $Tf$ by the formula
$$
Tf=\sum\limits_{k=0}^{n-1}\sum\limits^{n-1}_{l=0} f(x'_k,x'_l)
\int\int\limits_{\Delta_{kl}}
\frac{d\tau_1 d\tau_2}{((\tau_1-t_1)^2+(\tau_2-t_2)^2)^\lambda}+R_{nn}(f).
\eqno (6.8)
$$

Consider another cubature formula for calculating the integral  $Tf.$

Let $(t_1,t_2) \in \Delta_{ij}. $ By $\Delta_*$
denote the union of the square  $\Delta_{ij}$ and of those squares
$\Delta_{kl}$ which  have common points with the $\Delta_{ij}$.
Consider the formula
$$
Tf=f(x_i',x_j')\int\int\limits_{\Delta_*}
\frac{d\tau_1 d\tau_2}{((\tau_1-t_1)^2+(\tau_2-t_2)^2)^\lambda}+
$$
$$
+\sum\limits_{k=0}^{n-1}\sum\limits^{n-1}_{l=0}{}' f(x'_k,x'_l)
\int\int\limits_{\Delta_{kl}}
\frac{d\tau_1 d\tau_2}{((\tau_1-t_1)^2+(\tau_2-t_2)^2)^\lambda}+R_{nn}(f),
\eqno (6.9)
$$
where  $\sum\sum{}'$ means summation over the squares which do not
belong to $\Delta_*.$

{\bf Theorem 6.4.} {\it Among all cubature formulas  (2.1)
with $\rho_1=\rho_2=0$ and $n_1=n_2=n$,  formula (6.8), with the error
estimate  (6.15),
is optimal with respect to order.}

{\bf Remark.} {\it Similar statement holds for formula (6.9).}

{\bf Proof of Theorem 6.4.}
Let us estimate errors of formulas (6.8)  and (6.9).

The error of formula (6.8) can be estimated as follows:
$$
|R_{nn}(f)| \le \sum\limits_{k=0}^{n-1}\sum\limits^{n-1}_{l=0}{}'
\int\int\limits_{\Delta_{kl}}\frac{|f(\tau_1,\tau_2)-f(x'_i,x'_j)|}
{((\tau_1-t_1)^2+(\tau_2-t_2)^2)^\lambda}d\tau_1 d\tau_2 +
$$
$$
+\sum\limits_{k=0}^{n-1}\sum\limits^{n-1}_{l=0}{}''
\int\int\limits_{\Delta_{kl}}\frac{|f(\tau_1,\tau_2)-f(x'_k,x'_l)|}
{((\tau_1-t_1)^2+(\tau_2-t_2)^2)^\lambda}d\tau_1 d\tau_2 =r_1+r_2,
\eqno (6.10)
$$
where  $\sum\sum'$ means summation over  $k$ and  $l$ such that
the squares  $\Delta_{kl}$ belong to $\Delta_*$, and $\sum\sum''$ means
summation over the other squares.

Let us estimate $r_1$ and $r_2$:
$$
r_1 \le \frac{2}{n^\alpha}\int\int\limits_{\Delta_{*}}
\frac{d\tau_1d\tau_2}{((\tau_1-t_1)^2+(\tau_2-t_2)^2)^\lambda}
\le \frac{c}{n^{2-2\lambda+\alpha}}=o(n^{-\alpha}); \eqno (6.11)
$$
$$
r_2 \le \frac{4}{1+\alpha}\frac{1}{n^{2+\alpha}}
\sum\limits_{k=0}^{n-1}\sum\limits^{n-1}_{l=0}{}'' h(\Delta_{kl}).\eqno (6.12)
$$
Here  $h(\Delta_{kl})$ denotes the maximum value
of the function $((\tau_1-t_1)^2+(\tau_2-t_2)^2)^{-\lambda}$ in the square
$\Delta_{kl}$.

One has:
$$
\left|\int\int\limits_{\Delta_{kl}}
\left[\frac{1}{((\tau_1-t_1)^2+(\tau_2-t_2)^2)^\lambda}-h(\Delta_{kl})\right]
d\tau_1d\tau_2\right|=
$$
$$
=\left|\int\int\limits_{\Delta_{kl}}
\left[\frac{1}{((\tau_1-t_1)^2+(\tau_2-t_2)^2)^\lambda}-\frac{1}
{((x_k-t_1)^2+(x_l-t_2)^2)^\lambda}\right]d\tau_1d\tau_2\right| \le
$$
$$
\le\int\int\limits_{\Delta_{kl}}\left|\frac{2\lambda(x_k-t_1+q_1(\tau_1-x_k))(\tau_1-x_k))}
{((x_k-t_1+q_1(\tau_1-x_k))^2+(\tau_2-t_2)^2)^{\lambda+1}}d\tau_1d\tau_2\right|+
$$
$$
+\int\int\limits_{\Delta_{kl}}\left|\frac{2\lambda(x_l-t_2+q_2(\tau_2-x_l))(\tau_2-x_l))}
{((x_k-t_1)+(x_l-t_2+q_2(\tau_2-x_l))^2)^{\lambda}}d\tau_1d\tau_2\right|\le
$$
$$
\left.\le\int\int\limits_{\Delta_{kl}}\frac{2\lambda(\tau_1-x_k)}
{((x_k-t_1+q(\tau_1-x_k))^2+(\tau_2-t_2)^2)^{\lambda+1/2}}d\tau_1d\tau_2\right|+
$$
$$
\left.+\int\int\limits_{\Delta_{kl}}\frac{2\lambda(\tau_2-x_l))}
{((\tau_1-t_1)^2+(x_l-t_2+q_2(\tau_2-x_l))^2)^{\lambda+1/2}}
d\tau_1d\tau_2\right|\le
$$
$$
\le\frac{2^4\lambda}{n^3}\frac{n^{2\lambda+1}}{(k^2+l^2)^{\lambda+1/2}}=
\frac{2^4\lambda}{n^{2-2\lambda}}\frac{1}{(k^2+l^2)^{\lambda+1/2}},
$$
where it was assumed  that $k \ge i+1,$ and $l \ge j+1.$  Estimates for
the other combinations of $k$ and $l$ are similar. Thus:
$$
r_2 \le \frac{1}{(1+\alpha)}\frac{1}{n^\alpha}
\sum\limits_{k=0}^{n-1}\sum\limits^{n-1}_{l=0}{}''
\int\int\limits_{\Delta_{kl}}\frac{d\tau_1 d\tau_2}
{((\tau_1-t_1)^2+(\tau_2-t_2)^2)^\lambda}+
$$
$$
+\frac{1}{(1+\alpha)}\frac{2^3\lambda}{n^{2-2\lambda+\alpha}}
\sum\limits_{k=0}^{n-1}\sum\limits^{n-1}_{l=0}{}''
\frac{1}{(k^2+l^2)^{\lambda+1/2}}.
\eqno (6.13)
$$

Let us estimate the last term in the above inequality.

One has:
$$
\sum\limits_{k=0}^{n-1}\sum\limits^{n-1}_{l=0}{}''
\frac{1}{(k^2+l^2)^{\lambda+1/2}}\le
\sum\limits_{k=-[n/2]}^{n/2}\sum\limits^{n/2}_{l=-[n/2]}{}^*
\frac{1}{(k^2+l^2)^{\lambda+1/2}} \le
$$
$$
{\bf \le c}
\left\{
\begin{array}{ccc}
1, \quad \lambda>1/2\\
\log n, \quad \lambda=1/2\\
n^{1-2\lambda}, \quad \lambda<1/2\\
\end{array}
\right.
 \eqno (6.14)
$$
where  $\sum\sum^*$ means summation over $k$ and $l$,
$(k,l) \ne (0,0).$

In deriving (6.14) we have used the known result  ([17], Theorem 56)
which says that a number of points with integer-value coordinates,
situated
in the circle $x^2+y^2=r^2$, is equal to $\pi r^2+O(r).$

From inequalities (6.13) and (6.14) it follows that
$$
r_2 \le \frac{(1+o(1))}{(1+\alpha)}\frac{1}{n^\alpha}
\sum\limits_{k=0}^{n-1}\sum\limits^{n-1}_{l=0}{}''
\int\int\limits_{\Delta_{kl}}\frac{d\tau_1 d\tau_2}
{((\tau_1-t_1)^2+(\tau_2-t_2)^2)^\lambda}.
$$
This and (6.11) yield:
$$
R_{nn}[H_{\alpha\alpha}(D)] \leq \frac{1+o(1)}{(1+\alpha)n^\alpha}
\sup_{(t_1,t_2)\in D}\int\limits_{-1}^1\int\limits_{-1}^1
\frac{d\tau_1 d\tau_2}
{((\tau_1-t_1)^2+(\tau_2-t_2)^2)^\lambda}\le
$$
$$
\le \frac{1+o(1)}{(1+\alpha)n^\alpha}
\int\limits_{-1}^1\int\limits_{-1}^1
\frac{d\tau_1 d\tau_2}
{(\tau_1^2+\tau_2^2)^\lambda}.
\eqno (6.15)
$$

Theorem 6.4 follows from a comparison the estimates
of  $\zeta [H_{\alpha,\alpha}(D)]$ and $R_{nn} [H_{\alpha,\alpha}(D)].$
$\blacksquare$

Let us construct optimal with respect to order cubature formula for 
calculating
integrals $Tf$ on the classes $W^{rr}.$
In the derivation of formula  (5.9) the local
spline $\varphi_{n}(t_1,t_2)$, approximating the function
$\varphi (t_1,t_2)$ in the domain $[0,2\pi;0,2\pi],$
was constructed.
A spline $f_{nn}(t_1,t_2)$, approximating the function $f(t_1,t_2)$
in the domain $[-1,1]\times [-1,1]$, can be constructed analogously.
Calculate the integral $Tf$ by the formula
$$
Tf=\int\limits_{-1}^{1}\int\limits_{-1}^{1}
\frac{f_{nn}(\tau_1,\tau_2)d\tau_1d\tau_2}
{((\tau_1-t_1)^2+(\tau_2-t_2)^2 )^\lambda}+
R_{nn}(f). \eqno (6.16)
$$

{\bf Theorem 6.5.}
{\it Let $\Psi= W^{r,r}(1), r=1,2,\dots,$ and calculate
the integral $Tf$ by  formula (2.1)
with $\rho_1=\rho_2=r-1$, and $ n_1=n_2=n.$
Then cubature formula (6.16), which has the error
$$
R_{nn}(\Psi)\leq (1+o(1))\frac{2R_{r1}(1)}{(r+1)!(n-1+[R_{r1}(1)]^{1/r})^r}
\int\limits_{-1}^{1}\int\limits_{-1}^{1}
\frac{d\tau_1d\tau_2}
{(\tau_1^2+\tau_2^2 )^\lambda},
$$
is optimal with respect to order.
Here
$R_{rq}(t)$ is a polynomial of degree $r$, least deviating from zero
in $L_q([-1,1]).$}

As in the proof of the Theorem 5.6
one gets the following estimate
$$
R_{nn}(\Psi)\leq (1+o(1))\frac{2R_{r1}(1)}{(r+1)!(n-1+[R_{r1}(1)]^{1/r})^r}
\int\limits_{-1}^{1}\int\limits_{-1}^{1}
\frac{d\tau_1d\tau_2}
{(\tau_1^2+\tau_2^2 )^\lambda}.
$$

Comparing this estimate with the estimate of $\zeta_{nn}[W^{r,r}(1)]$
from Theorem 6.3 one finishes the proof.
$\blacksquare$

{\bf 7. Calculation of weakly singular integrals
               on piecewise continuous surfaces.}

In Sections 5 and 6 asymptotically optimal methods
for calculating weakly singular integrals defined on the squares
$[0,2\pi]^2$
or $[-1,1]^2$ were constructed.

It is of interest to study optimal methods for calculating weakly singular
integrals on piecewise-Lyapunov surfaces.

Consider the integral
$$
   Jf=\iint\limits_{G}\frac{f(\tau_1,\tau_2,\tau_3)dS}
   {\Bigl((\tau_1-t_1)^2+(\tau_2-t_2)^2+(\tau_3-t_3)^2\Bigr)^{\lambda}},
   \;t_1,t_2,t_3\in G,  \eqno (7.1)
$$
where G is a Lyapunov surface of class $L_s(B,\alpha).$

We show that the results derived in Sections 5 and 6 can be partially
generalized to the integrals (7.1).

Calculate integrals (7.1) by the formula:
$$
   Jf=\sum\limits_{k=1}^{n}\sum\limits_{|v|=0}^{\rho}
   p_{kv}f^{(v)} (M_k)+R_n(f,G,M_k,p_{kv},t),  \eqno (7.2)
$$
where $t=(t_1,t_2,t_3),\; v=(v_1,v_2,v_3),\; |v|=v_1+v_2+v_3,$
$ f^{(v)}(t_1,t_2,t_3)=\frac{\partial^{|v|}f}{\partial t_1^{v_1}
\partial t_2^{v_2}\partial t_3^{v_3}}.$

The error of formula (7.2) is:
$$ R_n(f,G,M_k,p_{kv})=\sup_{t\in G}|R_n(f,G,M_k,p_{kv},t)|.$$

Assume $f\in \Psi_1,$ and $G\in \Psi_2.$
Then the error of formula (7.2) on the classes
$\Psi_1$ and $\Psi_2$ is:
$$ R_n(\Psi_1,\Psi_2,M_k,P_{kv})=\sup_{f\in \Psi_1,G\in\Psi_2}
R_n(f,G,M_k,p_{kv}).$$

Let
$$
\zeta_n[\Psi_1,\Psi_2]:=\inf_{M_k,p_{kv}}R_n(\Psi_1,\Psi_2,M_k,p_{kv}).$$

A cubature formula with nodes $M_k^*$ and weights $p_{kv}^*$ is
called optimal, asymptotically optimal, optimal with respect to
order on the class of functions $\Psi_1$ and surfaces $\Psi_2,$ if
$$ \frac{R_n(\Psi_1,\Psi_2,M_k^*,p_{kv}^*)}
        {\zeta_n[\Psi_1,\Psi_2]}=1,\sim 1, \asymp 1,$$
respectively.

Let $\Psi_1=H_{\alpha}(1),\,\, 0<\alpha\leqslant 1,$ and
$\Psi_2=L_1(B,\beta)\,\, 0<\beta\leqslant 1.$
Let us construct an optimal
with respect to order method for calculating integrals
(7.1) on the classes of functions $\Psi_1$ and surfaces
$\Psi_2.$ Let $S(G)$ be a "square" of the surface $G.$ Divide the
surface $G$ into $n$ parts $g_k,\; k=1,2,\ldots,n,$ so that a
"square" of each of the domains  $g_k$ has the area of order $|S(G)|/n,$
where $|S(G)|$ is the area of $S(G)$.
We take a
point $M_k$ in each of domains $g_k$
at the center of the domain $g_k$.

Calculate integral (7.1) by the formula
$$
   Jf=\sum_{k=1}^{n}f(M_k)\iint\limits_{g_k}\frac{dS}
   {\Bigl((\tau_1-t_1)^2+(\tau_2-t_2)^2+(\tau_3-t_3)^2\Bigr)^{\lambda}}
   +R_n(f,G). \eqno (7.3)
$$
{\bf Theorem 7.1. } {\it Formula (7.3), has the
error $$R_n(\Psi_1,\Psi_2)\asymp n^{-\alpha/2},$$ and is optimal with
respect to order on the classes
$\Psi_1=H_{\alpha},\,\,0<\alpha\leqslant 1,$ and
$\Psi_2=L_1(B,\beta),\,\,\,0<\beta\leqslant 1,$ among all
formulas (7.2) with $\rho=0.$}

{\bf Proof. } Assume for simplicity that the surface $G$ is given by the
equation $z=\varphi(x,y),\;(x,y)\in G_0,$ $\varphi(x,y)\geqslant0$.
Let  $\varphi_x(x,y):=p,\,\, \varphi_y(x,y):=q$.
Write the integral $Jf$ as
$$
   Jf=\iint\limits_{G_0}\frac{f(\tau_1,\tau_2,\varphi(\tau_1,\tau_2))
   \sqrt{1+p^2(\tau_1,\tau_2)+q^2(\tau_1,\tau_2)}d\tau_1d\tau_2}
   {\Bigl[(\tau_1-t_1)^2+(\tau_2-t_2)^2+
   (\varphi(\tau_1,\tau_2)-\varphi(t_1,t_2))^2\Bigr]^{\lambda}}.
\eqno (7.4)
$$

The function $f(\tau_1,\tau_2,\varphi(\tau_1,\tau_2))$ belongs to
the H\"older class $H_{\alpha}$ over $G_{0},$ and the function
$\frac{\sqrt{1+p^2+q^2}}{\Bigl[(\tau_1-t_1)^2+(\tau_2-t_2)^2+
(\varphi(\tau_1,\tau_2)-\varphi(t_1,t_2))^2\Bigr]^{\lambda}}$ is
positive.

Let $M_k=(m_1^k,m_2^k,m_3^k)$ be the nodes of cubature formula
(7.2). Let
$\psi(\tau):=(d(\tau,\{M_k\}))^{\alpha},$ where $d(\tau,\{M_k\})$ is
the distance between the point $\tau$ and the set of the nodes
$\{M_k\},$ where  the distance is measured along the geodesics
 of the surface $G.$
This distance satisfies the H\"older condition
$H_{\alpha}(1).$
Hence the function
$\psi^*(\tau_1,\tau_2)=\psi(\tau_1,\tau_2,\varphi(\tau_1,\tau_2))$
belongs to the H\"older class $H_{\alpha}(A)$ and vanishes at the
nodes $(m_1^k,m_2^k),$ $k=1,2,\ldots,n.$ Thus,

$$ \zeta_n(\Psi_1,\Psi_2)\geqslant
\frac{1}{S(G_0)}\iint\limits_{G_0}\iint\limits_{G_0}
\frac{\psi(\tau_1,\tau_2,\varphi(\tau_1,\tau_2))
   \sqrt{1+p^2+q^2}d\tau_1d\tau_2dt_1 dt_2}
   {\Bigl[(\tau_1-t_1)^2+(\tau_2-t_2)^2+
   (\varphi(\tau_1,\tau_2)-\varphi(t_1,t_2))^2\Bigr]^{\lambda}}\geqslant$$
$$\geqslant\frac{1}{S(G_0)}\iint\limits_{G_0}\psi(\tau_1,\tau_2,
\phi(\tau_1,\tau_2))d\tau_1d\tau_2\times
$$
$$
\times \min_{t}\iint\limits_{G_0}\frac{\sqrt{1+p^2+q^2}}
{\Bigl[(\tau_1-t_1)^2+(\tau_2-t_2)^2+
   (\varphi(\tau_1,\tau_2)-\varphi(t_1,t_2))^2\Bigr]^{\lambda}}
   d\tau_1d\tau_2\geqslant$$
$$\geqslant\frac{A}{n^{\alpha/2}}\min\limits_t\iint\limits_{G}
\frac{ds}{\bigl(r(t,\tau)\bigr)^{\lambda}},$$ where $S(G_0)$ is the
"square" of the surface $G_0.$

Therefore the error of formula
(7.3) is estimated by the inequality
$R_n\leqslant\frac{A}{n^{\alpha/2}}.$

Theorem 7.1 is proved.
$\blacksquare$

{\bf Remark 1. } {\it The method of decomposition of the domain $G$
into smaller parts $g_k,\;k=1,2,\ldots,n,$  described below,
is optimal with respect to
order for classes of functions
$\Psi_1=H_{\alpha}, \,\,0<\alpha\leqslant 1$, and of surfaces
$\Psi_2=L_0(B,\beta),\,\,\,0<\beta\leqslant 1$ for
$\alpha\leqslant\beta.$}

{\bf Remark 2. } {\it From formula (7.4) it follows that if the
function $f\in W^{r,r}(1)$ and the surface $G\in L_s(B,\alpha),$
then the function $f(\tau_1,\tau_2,\varphi(\tau_1,\tau_2))\in
W^{v,v}(A),$ where $v=\min(r,s).$
Therefore, repeating the above arguments, one
proves that the accuracy of calculation of integral
(7.4) by cubature formulas using $n$ values of integrand
function does not exceed  $O(n^{-v/2}).$}

From this remark it follows that if the surface $G$
consists of several parts, for example of surfaces $G_1$ and $G_2$
having common edge $L,$ then it is necessary to calculate the integrals
for the surface $G_1$ and the surface $G_2$ separately. If the
surface $G$ is divided into smaller parts
$g_k,\;k=1,2,\ldots,n,$ the domains $g_k$, the curve $L$ passes
inside of these domains, should be associated with the class of
surfaces $L_0(B,1).$ In these domains the accuracy of calculation
of the integral does not exceed  than $O(n_k^{-1}),$ where $n_k$
is the number of nodes of the cubature formula used in the
domain $g_k.$

For this reason the cusps and the nodes, in which three or more
domains $G_k$, which are parts of the domain $G$ touch each other, must
belong to the boundaries of the covering domains
$g_k,\;k=1,2,\ldots,n.$

The universal code for computing
the capacitances, described in Section 9, is
based on optimal with respect to order cubature formulas for
calculating integrals on the classes of functions
$H_{\alpha}, \,\,\, 0<\alpha\leqslant 1$ and of surfaces
$L_0(B,\beta),\,\,\,B=const,\; \alpha\leqslant\beta,\;
\beta\leqslant 1.$

The algorithm constructed in Section 9 is optimal on
this class of surfaces and does not require  special treatment
of edges and conical points of the surface.

When one studies cubature formulas on the classes
$W^{r,r}(A),$ $r>1$, and
$L_s(B,\beta),$ $s\geqslant
1,\;0\leqslant\beta\leqslant1,$ one has to develop
a method to compute accurately the integrals in a neighborhood of the
above singular points of the surface.

{\bf 8. Calculation of weights of cubature formulas.}

In calculating weakly singular integrals by cubature formulas (6.8)
it is necessary to calculate integrals of the form of
$$
J_{kl}(t_1,t_2)=\int\limits_{\Delta_{kl}}\frac
{d\tau_1d\tau_2}{\left((\tau_1-t_1)^2+(\tau_2-t_2)^2\right)^
\lambda}
$$
for different values  $(t_1,t_2)\in [-1,1]^2.$

Let  $(t_1,t_2)\in \Delta_{ij}.$ Let us consider two possibilities:
1) the square $\Delta_{kl}$ and the square $\Delta_{ij}$ have nonempty
intersection;
2) the square $\Delta_{kl}$ is does not have common points with
the square $\Delta_{ij}$.

First consider the second case,
when the function
$$
\varphi(\tau_1,\tau_2)=\frac 1{\left((\tau_1-t_1)^2+(\tau_2-t_2)^2
\right)^\lambda},
$$
is smooth. Here $(\tau_1,\tau_2)\in \Delta_{kl},$ and $(t_1,t_2)\in
\Delta_{ij}.$

In this case one has:
$$
\left|\frac{\partial^r\varphi(\tau_1,\tau_2)}{\partial\tau_1^r}
\right|\le \frac{r!2^{2r}}{\left((\tau_1-t_1)^2+(\tau_2-t_2)^2\right)
^{\lambda+r/2}}
$$
and, if the squares  $\Delta_{kl}$ and $\Delta_{ij}$
do not have common points, one gets:
$$
\left|\frac{\partial^r\varphi(\tau_1,\tau_2)}{\partial\tau_1^r}
\right|\le\frac{2^rr!n^{2\lambda+r}}{2^\lambda}.
$$

Similar estimates holds for partial derivative with respect to $\tau_2.$

Calculate the integral $J_{kl}(t_1,t_2)$ by the Gauss cubature
formula:
$$
J_{kl}(t_1,t_2)=\int\limits_{\Delta_{kl}}P_{mm}\left[\frac
1{\left((\tau_1-t_1)^2+(\tau_2-t_2)^2\right)^\lambda}\right]
d\tau_1d\tau_2+R_{mm}(\Delta_{kl}),
$$
where $P_{mm}=P_m^{\tau_1}P_m^{\tau_2},\ P_m^{\tau_i}\ (i=1,2) $ is
the projection
operator onto the set of interpolation polynomials of degree $m$
with nodes at the zeros of the Legendre polynomial, which maps the
segment  $[-1,1]$
onto the segment $[x_k,x_{k+1}]$ for $i=1$, and onto the segment
$[x_l,x_{l+1}]$ for $i=2.$

An integer $m$ is chosen so that
$|R_{mm}|\leq n^{-2-\alpha}$ for cubature formulas
on the H\"older class $H_{\alpha\alpha}$, and $|R_{mm}|\leq n^{-r-\alpha}$
for cubature
formulas on the class $W^{rr}$.

This requirement is made because the error of calculation
of the coefficients $J_{kl}(t_1,t_2)$ must not exceed the error of
formula  (6.5).

Using $r$ derivatives of the integrand  in the error
$R_{mm}(\Delta_{kl})$, one gets:
$$
\left|R_{mm}(\Delta_{kl})\right|\le \frac{B_r2^rr!}{m^{r-1}}\left(\frac2n
\right)^{2-2\lambda},
$$
where  $B_r $ is the constant appearing in Jackson's theorems.
It is known that the constants  $B_r$ are bounded by a constant, denoted
$b$, uniformly with respect to $r$.
In the case of periodic functions $b=1$
([18]), and in the general case $b$ is apparently unknown.

If $r=2$ and
$m=B_r2^rr!n^{2\lambda}$, then one gets the error estimate given for cubature formula
(6.5).

Now, consider a method for calculating the integrals
$J_{kl}(t_1,t_2)$ when the square $\Delta_{kl}$ has nonempty
intersection with the square $\Delta_{ij}.$ For definiteness we consider
the calculation of the integral $J_{ij}(t_1,t_2)$
by the formula:
$$
J_{ij}(t_1,t_2)=\int\limits_{\Delta_{ij}}P_{mm}\left[\frac
1{\left((\tau_1-t_1)^2+(\tau_2-t_2)^2\right)^\lambda+h}\right]
d\tau_1d\tau_2+R_{mm}(\Delta_{ij}),
$$
where $h=const>0$ will be specified below.

One has:
$$
|R_{mm}(\Delta_{ij})|\le h\int\limits_{\Delta_{ij}}\frac
{d\tau_1d\tau_2}{\left((\tau_1-t_1)^2+(\tau_2-t_2)^2\right)^\lambda
\left(\left((\tau_1-t_1)^2+(\tau_2-t_2)^2\right)^\lambda+h\right)}+
$$
$$
+\int\limits_{\Delta_{ij}}D_{mm}\left[\frac 1{\left((\tau_1-t_1)^2+
(\tau_2-t_2)^2\right)^\lambda+h}\right]d\tau_1d\tau_2=r_1+r_2,
$$
where $D_{mm}=I-P_{mm},$ and $ I $ is an identity operator, and
$$
r_1\le h\int\limits_{\Delta_{ij}}\frac {d\tau_1d\tau_2}
{\left((\tau_1-t_1)^2+(\tau_2-t_2)^2\right)^{(\lambda+1)/2}
\left(\left((\tau_1-t_1)^2+(\tau_2-t_2)^2\right)^\lambda+h\right)^
{(1+\lambda)/2}}\le
$$
$$
\le\frac{2\pi}{1-\lambda}h^{(1-\lambda)/2}
\left(\frac{2\pi}n\right)^{1-2\lambda}. \eqno (8.1)
$$

The function $\frac 1{\left((\tau_1-t_1)^2+(\tau_2-t_2)^2\right)^\lambda+h}$
is infinitely smooth.
Using bounds for its first  derivatives for
$\lambda \ge 1/2$, one gets:
$$
r_2\le \frac{8\lambda B_1}{n^4h^2m}. \eqno (8.2)
$$

From inequality (8.1) it follows that for getting accuracy
$O(n^{-1-\alpha})$
one has to have  $h=n^{-2(2\lambda + \alpha)/(1-\lambda)}$
and from inequality (8.2) it follows that
one has to have  $m=max([n^{(8\lambda + 4\alpha)/(1-\lambda)+\alpha-3}],1)$.

{\bf 9. Iterative methods for calculating electrical capacitancies of
conductors of arbitrary shapes.}

Numerical methods for solving  electrostatic problems, in particular,
calculating capacitancies of conductors of arbitrary shapes,
are of practical interest in many applications.
Electrostatic problems solvable in closed form are collected in
[19,20,21].
Some of the problems were solved in closed form using  integral equations,
Wiener-Hopf  and singular integral equations [22]. Electrostatic problems
for a finite circular hollow cylinder (tube) were studied in [23] by numerical
methods. In [24] the variational methods of Ritz and Trefftz are
discussed.
Galerkin's and other projection methods are studied in [25].
In practice these methods are time-consuming and variational methods in
three-dimensional static problems probably have some advantages over
the grid method. There exists a vast literature on calculation of the
capacitances of perfect conductors [20,26]. In [20] there is a reference
section which gives the capacitance of the conductors of certain shapes
(more than 800 shapes are considered in [20]).
In [26] and [27] a systematic exposition of variational methods for
estimation
of the capacitances is given. In [28] there are some programs for calculating
the two-dimensional static fields using integral equations method.
In monograph [7]  iterative methods for solving interior and exterior
boundary value problems in electrostatics are proposed and mathematically
justified. Upper and lower estimates for some functionals of electrostatic
fields are obtained in [7]. Such functionals are the capacitances of 
perfect
conductors  and the
polarizability tensors of  bodies of arbitrary
shape. These bodies are described by their dielectric permittivity,
magnetic permeability and conductivity. They can be homogeneous or flaky.
The main point is: these bodies have arbitrary geometrical shapes.

The methods, developed in [7], yield analytical formulas for
calculation of the capacitances and polarizability tensors of
bodies of arbitrary shapes with any given accuracy. Error estimates for
these formulas are obtained in [7].
We give here the formulas for calculating the capacitances
of the conductors of arbitrary shapes [7]:

%\begin{equation}
 $$       C^{(n)} = 4 \pi \varepsilon_0 S^2
        \left\{ \frac{(-1)^n}{(2 \pi)^n} \int_\Gamma \int_\Gamma
        \frac{dsdt}{r_{st}}
     \underbrace{\int_\Gamma\dots
     \int_\Gamma}_{n \hbox{\begin{tiny}\ times\end{tiny}}}
        \psi(t,t_1) \dots \psi (t_{n-1}, t_n) dt_1 \dots dt_n
        \right\}^{-1}_,
$$
  %      \end{equation}
%\begin{equation}
where $S$ is the surface area of the surface $\Gamma$ of the conductor,
$\varepsilon_0$ is the dielectric constant of the medium,
$r_{st}:=|s-t|$, and $\psi(t,s):=\frac{\partial}{\partial N_t}
        \frac{1}{r_{st}},$

 $$       C^{(0)} = \frac{4 \pi \varepsilon_0 S^2}{J} \leq C, \quad
J \equiv \int_\Gamma \int_\Gamma \frac{dsdt}{r_{st}}, \quad
        S = \hbox{meas} \Gamma.
 $$

It is proved in [7] that
 $$       \left|C - C^{(n)}\right| \leq Aq^n, \quad 0 < q < 1,
   $$
where $A$ and $q$ are constants which depend only on the geometry of
$\Gamma$.

We use these formulas are used  to construct the
computer code for calculating the capacitances of the conductors of
arbitrary shapes.

It is proved in [7], that
$$
C^{(n)}=4\pi \varepsilon_{0} S^2 \left(\int\limits_{\Gamma}
\int\limits_{\Gamma}
r_{st}^{-1} \delta_n(t) dt ds\right)^{-1},   \eqno (9.1) $$
where $\delta_n$ is defined by the iterative process:
$$
\delta_{n+1}=-A\delta_n, \  \delta_0=1,  \int\limits_{\Gamma} \delta_n dt
=S,  \eqno(9.2) $$
and $A$ is defined by the formula:
$$
A\delta=\int\limits_{\Gamma} \delta(t) \frac{\partial}{\partial N_s}
\frac{1}{2\pi r_{st}} dt, $$
where $N_s$ is the outer unit normal to $\Gamma$ at the point $s.$

To use iterative process (9.2), one has to calculate the weakly
singular integral
$$
\frac{1}{2\pi} \int\limits_{\Gamma} \delta(t) \frac{\partial}{\partial N_S}
\frac{1}{r_{st}} dt. \eqno (9.3)$$
Let us  describe the construction of a cubature formula for
calculating integral (9.3),
assuming for simplicity that the domain $G$, bounded by the surface
$\Gamma$, is convex.
This asumption can be removed.

Let $\mathbb{S}$ be the inscribed in the conductor sphere of maximal
radius $r^*$, centered at the origin. Introduce the
spherical coordinates system $(r,
\phi,\theta),$   and the set of the nodes $(r^*,\phi_k, \theta_l),$
where $\phi_k=2 k \pi/n, \ k=0,1,\dots,n, \ \theta_l=\pi l/m, \ l=0,1,
\dots,m.$ Assume that $m$ is
even, and cover the sphere $\mathbb{S}$ with the spherical triangles
$\Delta_k, \
k=1,2,\dots,N, \ N=2n(m-1).$

Let us describe the construction of the spherical triangles.
For $0\le \Theta\le
\pi/m$ the triangles $\Delta_k, k=1,2,\dots,n$ have vertices  $(r^*,0,0), $
$(r^*,\phi_{k-1},\theta_1),$ $ (r^*,\phi_k,\theta_1),$ $ k=1,2,\dots,n.$

For $\theta_l \le \theta \le \theta_{l+1}, \ l=1,2,\dots, m/2-1,$ the triangles
$\Delta_k, k=n+2 n (l-1)+j, \ 1\le j \le 2n$ are constructed
as follows.
The rectangle $[0,2\pi;\theta_l,\theta_{l+1}]$ is covered with the
squares $\Delta_{kl}=[\phi_k,\phi_{k+1};\theta_l,\theta_{l+1}], \
k=0,1,\dots,n-1.$ Each of the squares $\Delta_{kl}$ is divived into two
equal triangles $\Delta_{kl}^1$ and $\Delta_{kl}^2, \ k=0,1,\dots,n-1,\
l=1,2,\dots,m/2-1.$ The spherical triangles
$\Delta_{kl}^1$ and $\Delta_{kl}^2, \ k=0,1,\dots,n-1,\
l=1,2,\dots,m/2-1,$ are images of  triangles  $\Delta_{kl}^1$ and
$\Delta_{kl}^2$ on the sphere $\mathbb{S}$

As a result of these constructions the sphere $\mathbb{S}$ is covered with
triangles $\Delta_k, \ k=1,2,\dots,N.$

We draw the straight lines through the origin and  vertices of
the triangle $\Delta_k, \ k=1,2,\dots,N.$
The points of intersection of these lines
with the surface $\Gamma$ are vertices of the triangle $\overline\Delta_k,
\ k=1,2,\dots,N.$ As a result of these constructions the surface
$\Gamma$ is approximated by the surface $\Gamma_N$  consisting of
triangle $\overline\Delta_k,\ k=1,2,\dots,N,$ and integral (9.3)
is approximated by the integral
$$
U(s)=\int\limits_{\Gamma_N} \delta(t) \frac{\partial}{\partial N_S}
\frac{1}{r_{st}} dt. \eqno (9.4)
$$

We fix each triangle $\overline\Delta_k,\ k=1,2,\dots,N,$ and
associate with it a point $\tau_k \in \overline\Delta_k,\ k=1,2,\dots,N,$
equidistant from the vertices of the triangle $\overline\Delta_k,\
k=1,2,\dots,N.$
We calculate integral (9.4) at the points $\tau_k, \ k=1,2,\dots,N,$
by the cubature formulas constructed in paragraphs 5-7
for the H\"older classes.
After calculating the values $U(\tau_k),\ k=1,2,\dots,N$ by these cubature
the integral
$$
\tilde C^{(1)}=-4\pi \varepsilon_{0} S_N^2 \left(\int\limits_{\Gamma_N}
\int\limits_{\Gamma_N} r_{st}^{-1} \tilde U(t) dt ds\right)^{-1}
$$
is calculated,
where $\tilde U(t)=U(\tau_k)$ for $t \in \overline\Delta_k,\ k=1,2,\dots,N,$
$S_N$ is area of the surface $\Gamma_N,$ $ \tilde C^{(1)}$ is
approximation to the value of $C^{(1)}.$
The successive iterations are calculated analogously.

{\bf 10. Numerical examples.}

In this section the numerical results are given.
As an example we calculated the capacitances of
various ellipsoids, because for ellipsoids one knows ([19])
the analytical formula for the capacitance, which makes it possible to
evaluate the accuracy of the numerical results.
Consider the ellipsoid:
$$
\frac{x^2}{a^2}+\frac{y^2}{b^2}+\frac{z^2}{c^2}=1.
$$
It is known  [19,20] that the exact value of the capacitance of
ellipsoid with $a=b$ is:
$$
C=\frac{4\pi\varepsilon_0 \sqrt{(a^2-c^2)}}{\arccos(c/a)}.$$

Let  $a=b=1,$ and $\varepsilon_0=1$. We have calculate the capacitance
$C$ for different values of the semiaxis $c.$ The results of
the calculations
are given in Table 1.

It is known ([7], p.43), that the capacitance of a metallic disc of
radius
$a$ is $C=8a\varepsilon_0$, and one can see from Table 1, that
asymptotically, as $c\to 0$, this formula can be used practically for
the
ellipsoids with $c\leq 0.001$ with the error about 0.005.

\newpage

\addtolength\oddsidemargin{-2cm}

%\size{11}{11pt}\selectfont%

\begin{center}
Table 1 \\
\end{center}

$$\vbox{\offinterlineskip
  \hrule
  \halign{\vrule#&\strut\quad\hfil#\hfil&\vrule#&\strut\quad\hfil#\hfil&\vrule#&\strut\quad\hfil#\hfil&
          \vrule#&\strut\quad\hfil#\hfil&\vrule#&\strut\quad\hfil#\hfil&\vrule#&\strut\quad\hfil#\hfil&
          \vrule#&\strut\quad\hfil#\hfil&\vrule#&\strut\quad\hfil#\hfil&\vrule#\cr
  height4pt&\omit&&\omit&&\omit&&\omit&&\omit&&\omit&&\omit&&\omit&\cr
  & \hfil C \quad && n \quad && m \quad && N \quad && Exact value \quad && Error \quad &&
   Relative error \quad && Calculation time \quad &\cr
  height4pt&\omit&&\omit&&\omit&&\omit&&\omit&&\omit&&\omit&&\omit&\cr
  \noalign{\hrule}
  height4pt&\omit&&\omit&&\omit&&\omit&&\omit&&\omit&&\omit&&\omit&\cr
  & 0.9 && 40 && 30 && 2320 && 12.144630 && -0.221200 && 0.018212 && 25 sec &\cr
  & 0.5 && 40 && 30 && 2320 && 10.392304 && -0.222042 && 0.021366 && 25 sec &\cr
  & 0.1 && 40 && 30 && 2320 && 8.5020638 && -0.301189 && 0.035425 && 25 sec &\cr
  & 0.01 && 40 && 30 && 2320 && 8.050854 && 0.072132 && 0.008959 && 25 sec &\cr
  & 0.001 && 40 && 30 && 2320 && 8.005092 && -0.821528 && 0.106374 && 25 sec &\cr
  & 0.0001 && 40 && 30 && 2320 && 8.000509 && -1.068178 && 0.133513 && 25 sec &\cr
  \noalign{\hrule}
  height4pt&\omit&&\omit&&\omit&&\omit&&\omit&&\omit&&\omit&&\omit&\cr
  & 0.9 && 50 && 40 && 3900 && 12.144630 && -0.180510 && 0.014801 && 1 min 15 sec &\cr
  & 0.5 && 50 && 40 && 3900 && 10.392304 && -0.185642 && 0.017860 && 1 min 15 sec &\cr
  & 0.1 && 50 && 40 && 3900 && 8.5020638 && -0.288628 && 0.033947 && 1 min 15 sec &\cr
  & 0.01 && 50 && 40 && 3900 && 8.050854 && -0.372047 && 0.046212 && 1 min 15 sec &\cr
  & 0.001 && 50 && 40 && 3900 && 8.005092 && -0.586733 && 0.073295 && 1 min 15 sec &\cr
  & 0.0001 && 50 && 40 && 3900 && 8.000509 && -0.933288 && 0.116653 && 1 min 15 sec &\cr
  \noalign{\hrule}
  height4pt&\omit&&\omit&&\omit&&\omit&&\omit&&\omit&&\omit&&\omit&\cr
  & 0.9 && 60 && 50 && 5880 && 12.144630 && -0.152009 && 0.012516 && 4 min &\cr
  & 0.5 && 60 && 50 && 5880 && 10.392304 && -0.160023 && 0.015391 && 4 min &\cr
  & 0.1 && 60 && 50 && 5880 && 8.5020638 && -0.283364 && 0.033328 && 4 min &\cr
  & 0.01 && 60 && 50 && 5880 && 8.050854 && 0.532250 && 0.061110 && 4 min &\cr
  & 0.001 && 60 && 50 && 5880 && 8.005092 && -0.391755 && 0.048939 && 4 min &\cr
  & 0.0001 && 60 && 50 && 5880 && 8.000509 && -0.880394 && 0.110042 && 4 min &\cr
  height4pt&\omit&&\omit&&\omit&&\omit&&\omit&&\omit&&\omit&&\omit&\cr}
  \hrule}$$

\newpage

\addtolength\oddsidemargin{+2cm}

{\bf References.}

{\bf 1. I.V. Boikov},{\it Optimal with respect to Accuracy Algorithms of
   Approximate Calculation of Singular Integrals,}
   Saratov State University Press, Saratov, 1983. (Russian).

{\bf 2. I.V. Boikov},{\it Passive and Adaptive Algorithms for the Approximate
   Calculation of Singular Integrals,} Ch. 1, Penza Technical State
   Univ. Press, Penza, 1995. (Russian).

{\bf 3. I.V. Boikov},{\it Passive and Adaptive Algorithms for the Approximate
   Calculation of Singular Integrals,} Ch. 2, Penza Technical State
   Univ. Press, Penza, 1995. (Russian).

{\bf 4. I.V. Boikov, N.F.Dobrunina and L.N.Domnin}, {\it Approximate Methods
    of Calculation of Hadamard Integrals and Solution of
    Hypersingular Integral Equations,} Penza Technical State Univ.
    Press, Penza, 1996. (Russian).

{\bf 5. I.V.  Boikov and T.I. Poljakova}, {\it Asymptotically optimal algorithms
for calculation Poisson integrals, Schwarz integrals and Cauchy type
integrals} Optimalnii metodi vichislenii i ix primenenie k obrabotke
informachii. Collection of works. Publishing house of Penza State
Technical University. Penza, 12 (1996), 17-48.
(Russian)

{\bf 6. A.G.Ramm}, {\it Electromagnetic wave scattering by small
bodies of arbitrary shapes,} in the book: ``Acoustic,
 electromagnetic and elastic scattering-Focus on
 T-matrix approach" Pergamon Press, N. Y. 1980.
 537-546. (editors V. Varadan and V.  Varadan). 

{\bf 7.  A.G. Ramm}, {\it  Iterative Methods for Calculating Static Fields and
    Wave Scattering by Small Bodies,} Springer - Verlag, New York,
     1982.

{\bf 8. N.S. Bakhvalov}, {\it Properties of Optimal Methods of Solution of
   Problem of Mathematical Physics,} Zhurn. Vych. Mat. i Mat. Fiz.,
   10, N3, (1970), 555-588.

{\bf 9.} {\it Theoretical Bases and Construction of Numerical Algorithms for
    The Problems of Mathematical Physics,} (Ed. K.I. Babenko), Nauka,
    Moscow,  1979. (Russian).

{\bf 10. I.F. Traub and H.Wozniakowski}, {\it A General Theory of Optimal
    Algorithms,}  Academic Press, N. Y., 1980.

{\bf 11. S.M.Nikolskii}, {\it Quadrature Rules,}  Nauka, Moscow,
1979. (Russian).

{\bf 12. G.G. Lorentz}, {\it Approximation of function,} Chelsea
Publishing Company, New York, 1986.

{\bf 13. N.M. Gunter}, {\it Theory of potential and its application to basic
problems of mathematical physics,} GITTL, Moscow, 1953. (Russian).

{\bf 14. P.J. Davis and P. Rabinovitz}, {\it Methods of numerical
integration}
2nd ed., Academic Press, New York, 1984.

{\bf 15. V.I. Krylov}, {\it Approximate calculation of integrals,}
Nauka, Moscow, 1967. (Russian)

{\bf 16. N.S. Bakhvalov}, {\it Optimal linear methods of approximation
operators on convex classes of functions,} Zh. Vychisl. Mat. i Mat. Fiz.,
   11, N4, (1971), 1014-1018.

{\bf 17. A.A. Bukhshtab},  {\it Theory of numbers,}
Gos.
uchebno-pedag. izd. Minister. prosveshenija of Russia, Moscow, 1960.

{\bf 18. N.P. Kornejchuk}, {\it Exact Conctants,} Nauka, Moscow, 1990.
(Russian).

{\bf 19. L. Landau and E. Lifschitz}, {\it  Electrodynamics of
Continuous Media,} Pergamon Press, N.Y., 1960.

{\bf 20. Ju. Jossel, E. Kochanov, and M. Strunskij}, {\it Calculation of
Electrical Capacitance,} Energija, Leningrad, 1963. (Russian).

{\bf 21. P. Morse  and M. Feshbach}, {\it Methods of Theoretical Physics,}
vols 1 and 2, Me Graw-Hill, N.Y., 1953.

{\bf 22. B. Noble}, {\it  Wiener-Hopf Methods for Solution of Partial Differential
Equations,}  Pergamon Press, N.Y., 1958.

{\bf 23. L. Wainstein}, {\it  Static Problems for Circular Hollow
Cylinder of Finite
Length,} J. Tech. Phys, 32 (1962) 1162-1173; 37 (1967), 1181-1188.

{\bf 24. N. Miroljubov}, {\it  Methods of Calculating of
Electrostatic fields,}
High School, Moscow, 1963. (Russian).

{\bf 25. M. Krasnoselskij, G. Vainikko, P. Zabreiko, Ja. Rutickij
and V. Stecenko}, {\it  Approximate Solution of Nonlinear Equations,}
Wolters-Noordhoff, Groningen, 1972 (1969).

{\bf 26. G. Polya and G. Szego}, {\it Isoperemetrical Inequalities in
Mathematical Physics,} Princeton Univ. Press, Princeton, 1951.

{\bf 27. L. Payne}, {\it Isoperemetrical inequalities and  thear application,}
SIAM Rev, 9, (1967), 453-488.

{\bf 28. O. Tosoni}, {\it Calculation of Electromagnetic Fields on Computers,}
Technika, Kiev, 1967. (Russian).

\end{document}